\documentclass{amsart}
\usepackage{amsmath}
\usepackage{graphicx}
\usepackage{amssymb}
\usepackage{epstopdf}
\usepackage{tikz-cd}
\usepackage{xcolor}
\usepackage{hyperref}
\usepackage{tikz}
\usepackage{tikz-cd} 
\usetikzlibrary{cd}
\usepackage{framed}
\usepackage[T1]{fontenc}
\usepackage{bm}

\newtheorem{theorem}{Theorem}[section]
\newtheorem*{thm*}{Theorem A}
\newtheorem*{thm*B}{Theorem B}
\newtheorem*{thm*C}{Theorem C}
\newtheorem{lemma}[theorem]{Lemma}
\newtheorem{proposition}[theorem]{Proposition}
\newtheorem{definition}[theorem]{Definition}

\title[Cyclotomic Hecke fields and $L$-values of $\mathrm{GL}(2)$ cusp forms]{Generation of cyclotomic Hecke fields by $L$-values of cusp forms on $\mathrm{GL}(2)$ with certain $\mathbb{Z}_p$ twist}
\author{Jaesung Kwon}
\email{jaesungkwon@snu.ac.kr}
\address{Seoul National University, 1, Gwanak-ro, Gwanak-gu, Seoul, Republic of Korea, 08826}
                         
\date{\today}

\begin{document}

\begin{abstract}
Let $F$ be a number field, $f$ an algebraic automorphic newform on $\mathrm{GL}(2)$ over $F$, $p$ an odd prime does not divide the class number of $F$ and the level of $f$.
We prove that $f$ is determined by its $L$-values twisted by Galois characters $\phi$ of certain $\mathbb{Z}_p$-extension of $F$. Furthermore, if $F$ is totally real or CM, then under some mild assumption on $f$, the compositum of the Hecke field of $f$ and the cyclotomic field $\mathbb{Q}(\phi)$ is generated by the algebraic $L$-values of $f$ twisted by Galois characters $\phi$ of certain $\mathbb{Z}_p$-extension of $F$.
\end{abstract}


\maketitle
\numberwithin{equation}{section}
\setcounter{tocdepth}{2}
\tableofcontents

\section{Introduction} 
Arithmetic of modular $L$-values is one of the central area in number theory, as it is interesting itself but also important arithmetic consequences arise.
Let $p$ be a rational prime.
Rohrlich \cite{rohrlich1984onl} proved the non-vanishing of modular $L$-values twisted by cyclotomic characters, namely, Galois characters of the cyclotomic $\mathbb{Z}_p$-extension of $\mathbb{Q}$. This result deduces that it is finitely generated that the Mordell-Weil group of a CM elliptic curve over the cyclotomic $\mathbb{Z}_p$-extension of $\mathbb{Q}$.
Luo-Ramakrishnan \cite{luo1997determination} proved that for a newform $f$ of level $\Gamma_0(N)$, the field generated by $\mathbb{Q}(\mu_{p^\infty})$ and the $L$-values twisted by the cyclotomic characters coincides the cyclotomic Hecke field $\mathbb{Q}_f(\mu_{p^\infty})$ of $f$, where $\mu_{p^\infty}$ is the set of all $p$-power roots of unity. Recently, Sun \cite{sun2018generation} shows that for a newform $f$ of level $\Gamma_0(N)$ and a tame character $\eta=\eta_p\eta^{(p)}$ such that $f\otimes\eta_p$ has no $p$-power inner twist, the $L$-value of $f$ twisted by $\eta\chi$ generates the Hecke field $\mathbb{Q}_g(\eta\chi)$ over $\mathbb{Q}(\eta^{(p)})$ for almost all cyclotomic characters $\chi$, where $\eta_p$ and $g$ are the $p$-part of $\eta$ and the $p$-new part of $f$, respectively. This result applies to the Hecke fields of modular forms in the Hida family.

As a generalization of the case of modular forms, arithmetic of $L$-values of automorphic forms on $\mathrm{GL}(2)$ with various character twists also has been studied broadly. 
Let $F$ be a number field and $f$ a Hecke eigenform on $\mathrm{GL}(2)$ over $F$.
Rohrlich \cite{Rohrlich:1989} proved that $L$-values of $f$ are non-vanishing for infinitely many Hecke character twists.
Friedberg-Hoffstein \cite{friedberg1995nonvanishing} proved that $L$-values of $f$ are non-vanishing for infinitely many quadratic twists.
Kwon-Sun \cite{kwon2020non} proved that $L$-values of $f$ are non-vanishing for almost all character twists $\phi$ in $\Xi_\mathfrak{p}$, where $\Xi_\mathfrak{p}$ is the set of $p$-power order primitive characters of Galois groups $\mathrm{Gal}(F_m/F)$, $F_m$ is the ray class field of conductor $\mathfrak{p}^m$, and $\mathfrak{p}$ is a prime ideal of $F$ lying above $p$ such that $\mathfrak{p}$ is coprime to the class number of $F$ and the level of $f$, the residue degree and ramification degree are $1$, and the Galois group $\mathrm{Gal}(F_\infty/F):=\varprojlim_{m}\mathrm{Gal}(F_m/F)$ has $\mathbb{Z}_p$-rank $1$ (we recall this in Definition \ref{coates:wiles:extension}), so that $F_\infty^\Delta$ is a $\mathbb{Z}_p$-extension of $F$ and $\phi\in\Xi_\mathfrak{p}$ is a character of $\mathrm{Gal}(F_\infty^\Delta/F)\cong\mathbb{Z}_p$, where $\Delta$ is the torsion subgroup of $\mathrm{Gal}(F_\infty/F)$. 
Kwon \cite{kwon2023bianchi} studies the $\mathfrak{p}$-adic $L$-function of $f$ for characters $\phi$ in $\Xi_\mathfrak{p}$ when $F$ is imaginary quadratic.

Note that the Hecke field $\mathbb{Q}_f$ of $f$ is a number field and there is an algebraic $L$-value $L_f(\xi)\in\mathbb{Q}_f(\xi)$ for finite order Hecke characters $\xi$ over $F$ (see Hida \cite{Hida:1994} and Section \ref{mainresult:section}). From this fact and the results Kwon-Sun \cite{kwon2020non} and Sun \cite{sun2018generation}, we can naturally guess that $L_f(\phi)$ generates the cyclotomic Hecke field $\mathbb{Q}_f(\phi)$ for the characters $\phi\in\Xi_\mathfrak{p}$.

The main result of this paper is that if $F$ is totally real or CM, then the algebraic $L$-value $L_f(\phi)$ generates the cyclotomic Hecke field $\mathbb{Q}_f(\phi)$ for almost all $\phi\in\Xi_\mathfrak{p}$, which is a partial generalization of Sun \cite{sun2018generation}.

\subsection{Main results}
Let us follow the notations of the previous paragraphs.
Then, the main results of the present paper are following:
\begin{thm*}[Non-vanishing of $L$-values, Theorem \ref{Galois:average:L:value:prop}]
The $L$-value of $f$ is non-vanishing for almost all character twists $\phi$ in $\Xi_\mathfrak{p}$.
\end{thm*}

Note that Theorem A is equal to Kwon-Sun \cite[Corollary 7.1]{kwon2020non}. We have fixed some errors of Kwon-Sun \cite{kwon2020non} and reprove this result in the present paper.

\begin{thm*B}[Determination of newforms, Theorem \ref{determination:Lvalue}]
Let $f_1$ and $f_2$ be newforms on $\mathrm{GL}(2)$ over $F$. If the algebraic $L$-values $L_{f_1}(\phi)$ and $L_{f_2}(\phi)$ coincide for infinitely many $\phi\in\Xi_\mathfrak{p}$ up to a constant which does not depend on $\phi$, then $f_1=f_2$.
\end{thm*B}

Furthermore, if $F$ is totally real or CM, then we obtain our main result:
\begin{thm*C}[Generation of cyclotomic Hecke fields, Theorem \ref{themaintheorem}]
Let $f$ be a newform on $\mathrm{GL}(2)$ over $F$ whose central character is trivial at the finite places and the Hecke field $\mathbb{Q}_f$ is totally real (Note that $\mathbb{Q}_f$ is either totally real or CM).
Let us assume that $F$ is totally real or CM.
Then, there is an element $\phi_0\in\Xi_\mathfrak{p}$ such that
$$
\mathbb{Q}_f\big(\phi_0,A_f(\phi,\phi_0)\big)=\mathbb{Q}_f(\phi)
$$
for almost all $\phi\in\Xi_\mathfrak{p}$, where $A_f(\phi,\phi_0):=L_f(\phi)/L_f(\phi_0)$, which is well-defined by Theorem A.
\end{thm*C}
Let us make a remark that the assumption on $F$ is necessary because the proof of Theorem C needs the existence of the Galois representation attached to $f$, which is guaranteed by Harris-Lan-Taylor-Thorne \cite{Harris2014on} in this case.

\subsection{Outline and key ideas} 
In Section \ref{review:automorphicform:section}, we summarize the theory of automorphic forms on $\mathrm{GL}(2)$ over number fields decribed in Miyake \cite{miyake1971onautomorphic} and Hida \cite{Hida:1994}. Also we give the summarized version of the proof of functional equations of automorphic $L$-functions (the equation (\ref{approx:func:eq:'}), Proposition \ref{Lftn:analytic:conti} and \ref{fricke:twist:commute}) in Kwon-Sun \cite{kwon2020non}. These  functional equations plays an important role to study the automorphic $L$-values.
 
In Section \ref{subsec:hecke:char}, we define the additive twists of Galois averages of Hecke characters and automorphic $L$-values (Definition \ref{galois:averages} and \ref{galois:average:Lvalue}), and compute estimations of these averages by using the {\it approximate functional equation} obtained in Section \ref{review:automorphicform:section}. 
One of the important inputs is {\it coherent cone decomposition} by Rorhlich \cite{Rohrlich:1989}, which gives a lower bound of the absolute value of the norm of non-zero elements in a residue class and gives an upper bound of the number of elements in the intersection of a residue class and the fundamental unit cone.
From this, we can extract each Fourier-Whittaker coefficients of $f$ from its twisted $L$-values (Proposition \ref{main:thm:1}).

Finally, in Section \ref{mainresult:section}, we obtain the main results. We reprove the non-vanishing result of Kwon-Sun \cite{kwon2020non} (Theorem \ref{Galois:average:L:value:prop}) and show the determination of newforms by its $L$-values (Theorem \ref{determination:Lvalue}). By using these results, we can prove the generation of cyclotomic Hecke fields by the $L$-values (Theorem \ref{themaintheorem}). The key idea comes from Sun \cite{sun2018generation}, who utilizes the properties of the Galois representations attached to modular forms. 

Let us give a brief proof of our main theorem. 
Let $f$ be a newform on $\mathrm{GL}(2)$ over $F$.
Let us denote by $L_\phi$ the field generated by algebraic $L$-values $L_f(\phi)$ of $f$ for $\phi\in\Xi_\mathfrak{p}$.
If we assume that the field $L_\phi$ does not contain the cyclotomic field $\mathbb{Q}(\phi)$ for infinitely many $\phi\in\Xi_\mathfrak{p}$, then the twisted averages $\mathrm{Tr}\big(\beta L_f(\phi)\big)$ of algebraic $L$-values are vanishing. 
For a moment, let us recall that we can extract each Fourier-Whittaker coefficients of $f$ from its algebraic $L$-values $L_f(\phi)$.
Therefore, we observe that a linear combination $\sum_{\tau,\varphi}c_{\tau,\varphi}f^\tau\otimes\varphi$ of twisted Hecke eigenforms are vanishing, where $\tau$ and $\varphi$ run through the $\mathbb{Q}$-embeddings of the Hecke field $\mathbb{Q}_f$ to $\mathbb{C}$ and the dual of $\mathrm{Gal}(F_{n_0}^\Delta/F)$ for some $n_0>0$, respectively. By using the theory of Galois representation, we deduce from $\sum_{\tau,\varphi}c_{\tau,\varphi}f^\tau\otimes\varphi=0$ that some distinct Galois conjugates of newforms $f^{\tau_1}$ and $f^{\tau_2}$ are equivalent (Lemma \ref{modq:galois:repn}), which is a contradiction due to the newform theory. In conclude, $L_\phi$ contains $\mathbb{Q}(\phi)$ for almost all $\phi\in\Xi_\mathfrak{p}$. On the other hand, by using the determination result, we can also show that $L_\phi$ contains the Hecke field $\mathbb{Q}_f$. So we are done.

\subsection{Notations}\label{subsection:notation}
Let us provide notations which appear frequently in this paper:
\begin{itemize}
\item
For an integer $M>0$, let us denote by $\mu_M$ the set of $M$-th roots of unity.
\item
Let $F$ be a number field. 
\item
Let $\mathbb{A}_F$ be the adele ring of $F$, 
$\mathfrak{d}_F$ the different ideal of $F$, $D_F$ the discriminant of $F$, 
$F_\mathbb{R}$ the ring of infinite adeles of $F$, $h_F$ the class number of $F$, $I_F$ the set of $\mathbb{Q}$-embeddings of $F$ into $\mathbb{C}$, $N$ the absolute norm map of $F$, and $O_F$ the integer ring of $F$.
\item
For a place $v$ of $F$, let us denote by $F_v$ the $v$-adic completion of $F$, $a_v$ the image of $a\in\mathbb{A}_F$ under the projection $\mathbb{A}_F\to F_v$.
\item
For a finite place $v$ of $F$, let us denote by $\mathrm{ord}_v$ the $v$-adic valuation map of $K_v^\times$, $O_v$ the valuation ring of $\mathrm{ord}_v$.
\item
For a finite place $v$ of $F$, let us fix a uniformizer $\varpi_v$ of $K_v^\times$.
\item 
For a finite place $v$ and a subset $A$ of $F$, let us write $\alpha\in A\cap O_v^\times$ if $\alpha\in A$ and $\alpha$ is coprime to $v$.
\item
For complex valued functions $g$ and $h$, let us write $g\ll_{\mathcal{D}} h$ if $|g|\leq c_\mathcal{D}|h|$ for some constant $c_\mathcal{D}>0$, which depends only on the datum $\mathcal{D}$.
\item
Let us denote by $\mathbf{e}(z):=e^{2\pi iz}$ for $z\in\mathbb{C}$. Let $\mathbf{e}_F$ be the character of $F\backslash\mathbb{A}_F$ characterized by $\mathbf{e}_F(z)=\mathbf{e}\big(\mathrm{Tr}_{F/\mathbb{Q}}(z)\big)$ for $z\in F_\mathbb{R}$.
\end{itemize}

\section{Cusp forms on GL(2) and its {\it L}-functions}\label{review:automorphicform:section}
In this section, we review the theory of cusp forms on $\mathrm{GL}(2)$ over $F$ and its $L$-functions.
All the settings and the results in this section come from Hida \cite{Hida:1994}, Kwon-Sun \cite{kwon2020non}, and Miyake \cite{miyake1971onautomorphic}.

\subsection{Cusp forms on GL(2) }\label{cuspform}
Let $\Sigma(\mathbb{R})$ and $\Sigma(\mathbb{C})$ be the set of real places and the set of complex places of $F$, respectively. Let us denote by ${\rm c}$ the complex conjugation. Let $\mathbf{k}=\sum_{\tau\in I_F}k_\tau\tau$ and $\mathbf{m}=\sum_{\tau\in I_F}m_\tau\tau$ be elements of $\mathbb{Z}[I_F]$ satisfying the following conditions:
\begin{itemize}
\item[(1)] $k_\tau\geq 2$, $k_\tau=k_{{\rm c}\tau}$, and $m_\tau\geq 0 \text{ for any }\tau\in I_F$.
\item[(2)] $k_\tau+2m_\tau$ does not depend on $\tau\in I_F$.
\end{itemize}
Let us denote by $L_{\mathbf{k}}:=\mathrm{Sym}^{2k_\tau-2}(\mathbb{C}^2)^{\otimes\Sigma(\mathbb{C})}$ and
$
C_{F,\infty}^+:=\mathrm{SO}_2(\mathbb{R})^{\Sigma(\mathbb{R})}\times\mathrm{SU}_2(\mathbb{C})^{\Sigma(\mathbb{C})},
$
which acts on $L_{\mathbf{k}}$.
Let $\mathfrak{N}$ be a non-zero integral ideal of $F$. 
Let us define an open subgroup $U_0(\mathfrak{N})$ of $\mathrm{GL}_2(\widehat{O}_F)$ by
$
U_0(\mathfrak{N}):=\{\left(\begin{smallmatrix} a & b \\ c & d \end{smallmatrix}\right)\in \mathrm{GL}_2(\widehat{O}_F) \mid c\in\widehat{\mathfrak{N}}\}$,
where $\widehat{A}:=A\otimes_\mathbb{Z}\widehat{\mathbb{Z}}$ for a $\mathbb{Z}$-subalgebra $A$ of $F$.
Let us denote by $k:=k_\tau+2m_\tau$.
Let $\chi$ be a Hecke character modulo $\mathfrak{N}$ of infinity type $2\mathbf{t}-\mathbf{k}-2\mathbf{m}=(2-k)\mathbf{t}$, where $\mathbf{t}:=\sum_{\tau\in I_F}\tau\in\mathbb{Z}[I_F]$.
Let $J$ be a subset of $\Sigma(\mathbb{R})$.

\begin{definition}\label{adel:cuspform:defn} 
A cusp form on $\operatorname{GL}_2(\mathbb{A}_F)$ of weight $(\mathbf{k},\mathbf{m})$, type $J$, level $\mathfrak{N}$, and central character $\chi$ is a $C^\infty$-function $f:\operatorname{GL}_2(\mathbb{A}_F)\rightarrow L_{\mathbf{k}}$ such that
\begin{itemize}
\item[(1)] 
$D_\tau f=\frac{1}{2}k_\tau(k_\tau-2)f$ for the Casimir operators $\{D_\tau\mid\tau\in I_F\}$ of $\mathrm{SL}_2$.
\item[(2)] For elements $\gamma\in\mathrm{GL}_2(F)$, $z\in \mathbb{A}_F^\times$, $g\in\mathrm{GL}_2(\mathbb{A}_F)$, $u=\left(\begin{smallmatrix} a & b \\ c & d \end{smallmatrix}\right)\in U_0(\mathfrak{N})$, and $c\in C_{F,\infty}^+$, $f$ satisfies the following:
$$
f(\gamma z g uc)(\mathbf{x})=\chi(z)\chi\Big(\prod_{v|\mathfrak{N}}d_v\Big)\mathbf{e}\Big( \sum_{\tau\in J}k_\tau\theta_\tau-\sum_{\tau\in \Sigma(\mathbb{R})-J}k_\tau\theta_\tau\Big)f(g)(c_\tau\mathbf{x}_\tau)_{\tau\in\Sigma(\mathbb{C})},
$$ 
where 
$c_\tau=\big(\begin{smallmatrix}
\mathrm{cos}(2\pi\theta_\tau) & -\mathrm{sin}(2\pi\theta_\tau) \\ \mathrm{sin}(2\pi\theta_\tau) & \mathrm{cos}(2\pi\theta_\tau)\end{smallmatrix}\big)$ for $\tau\in\Sigma(\mathbb{R})$.
\item[(4)] Let $dx$ be a Haar measure on $F\backslash\mathbb{A}_F$. For $g\in \mathrm{GL}_2(\mathbb{A}_F)$,
$$
\int_{F\backslash \mathbb{A}_F}f\big(\left(\begin{smallmatrix} 1 & x \\ 0 & 1 \end{smallmatrix}\right) g\big)dx=0.
$$
\end{itemize}

From now on, let us denote by $S_{(\mathbf{k},\mathbf{m}),J}(\mathfrak{N},\chi)$ the space of such cusp forms $f$.
\end{definition}


\subsection{Hecke eigenforms and newforms on GL(2)}\label{newform}
Cusp forms $f$ in the space $S_{(\mathbf{k},\mathbf{m}),J}(\mathfrak{N},\chi)$ are related to harmonic differential forms on some $C^\infty$-hyperbolic manifold $Y$ so that there is a Fourier-Whittaker expansion of $f$:
\begin{definition}
For $f\in S_{(\mathbf{k},\mathbf{m}),J}(\mathfrak{N},\chi)$ and an ideal $\mathfrak{a}$ of F, let us denote by $a_f(\mathfrak{a})$ the Fourier-Whittaker coefficient of $f$ at $\mathfrak{a}$.
\end{definition}
A cusp form $f$ on $\mathrm{GL}(2)$ is called {\it Hecke eigenform} if $f$ is a common eigenfunction of all the Hecke operators.
Note that $a_f(O_F)\neq 0$ if $f$ is a Hecke eigenform (Miyake \cite[Lemma 2]{miyake1971onautomorphic}).
A Hecke eigenform $f$ is said to be {\it normalized} if $a_f(O_F)=1$. 

\begin{definition}
For normalized Hecke eigenforms $f_1\in S_{(\mathbf{k},\mathbf{m}),J}(\mathfrak{N},\chi_1)$ and $f_2\in S_{(\mathbf{k},\mathbf{m}),J}(\mathfrak{N},\chi_2)$ with the eigensystems $\lambda_1,\lambda_2$,
let us denote by $f_1\sim f_2$ if 
$$
\lambda_1(v)=\lambda_2(v)
$$ 
for almost all places $v$ of $F$. 
\end{definition}

The space $S_{(\mathbf{k},\mathbf{m}),J}(\mathfrak{N},\chi)$ is a $\mathbb{C}$-vector space equipped with an inner product (see Miyake \cite[(1.2)]{miyake1971onautomorphic}).
Let us denote by
$S^{\rm new}_{(\mathbf{k},\mathbf{m}),J}(\mathfrak{N},\chi)$ the orthogonal complement of the oldspace $S^{\rm old}_{(\mathbf{k},\mathbf{m}),J}(\mathfrak{N},\chi)$ (see Miyake \cite{miyake1971onautomorphic} for the definition of the oldspace).
\begin{definition}
$f\in S_{(\mathbf{k},\mathbf{m}),J}(\mathfrak{N},\chi)$ is called {\it newform} if $f$ is a normalized Hecke eigenform and $f\in S^{\rm new}_{(\mathbf{k},\mathbf{m}),J}(\mathfrak{M},\chi^\prime)$ for some divisor $\mathfrak{M}$ of $\mathfrak{N}$ and central character $\chi^\prime$ of modulus $\mathfrak{M}$ which induces $\chi$. 
\end{definition}
Note that if $f$ is a newform with the eigensystem $\lambda_f$, then $a_f=\lambda_f$.
If $f_1$ and $f_2$ are newforms such that $f_1\sim f_2$, then $f_1=f_2$ by the strong multiplicity one theorem (\cite[Corollary 3]{miyake1971onautomorphic}). 
\begin{lemma}\label{newform:lemma}
Let $f_i\in S_{(\mathbf{k},\mathbf{m}),J}(\mathfrak{N},\chi_i)$ $(i=1,\cdots,m)$ be a normalized Hecke eigenform.
If $f_i\not\sim f_j$ for each $i\neq j$, then $f_1,\cdots,f_m$ are $\mathbb{C}$-linearly independent.
\end{lemma}

\begin{proof}
Let us assume the contrary, then there is a minimal non-trivial relation $\sum_{i\in I}c_if_i=0$. Note that $|I|>1$ as $f_i\neq 0$ for each $i$. Choose distinct two elements $j,l\in I$. Since $f_j\not\sim f_l$, there is a place $v$ of $F$ such that $\lambda_j(v)\neq\lambda_l(v)$, where $\lambda_i$ is the eigensystem of $f_i$. Hence, we obtain that $\sum_{i\in I-\{j\}}c_i\big(\lambda_i(v)-\lambda_j(v)\big)f_i=0$, however $c_l\big(\lambda_l(v)-\lambda_j(v)\big)\neq 0$, which is contradict to the minimality of $I$.
\end{proof}

\subsection{{\it L}-function of cusp forms and its integral representation} 
From now on, $f$ is an element of $S_{(\mathbf{k},\mathbf{m}),J}(\mathfrak{N},\chi)$.

\begin{definition}[Automorphic $L$-function]
Let us define $L$-function of $f$ as follows:
$$
L(s,f):=\sum_{\mathfrak{a}}\frac{a_f(\mathfrak{a})}{N(\mathfrak{a})^s},\ \mathrm{Re}(s)>\frac{k}{2}+1,
$$
where $\mathfrak{a}$ runs through the ideals of $F$. Note that $a_f(\mathfrak{a})=0$ unless $\mathfrak{a}$ is integral by Hida \cite[Theorem 6.1]{Hida:1994}.
\end{definition}

Let $\widetilde{f}$ be the $\prod_{\tau\in\Sigma(\mathbb{C})}\binom{2k_\tau-2}{k_\tau-1} X_\tau^{k_\tau-1}Y_\tau^{k_\tau-1}$-component of $f$.
Let $F_{\mathbb{R},+}^\times$ be the set of totally positive elements of $F_{\mathbb{R}}^\times$. Let us denote by $E_F:=O_F^\times\cap F_{\mathbb{R},+}^\times$.
For a finite idele $a$ of $F$ and an element $y\in E_F\backslash F_{\mathbb{R},+}^\times$, let us define
$$
f_{a}(y):=|a|^{-1}|y|^{-\frac{k}{2}}\widetilde{f}\big(\left(\begin{smallmatrix} a & 0 \\ 0 & 1 \end{smallmatrix}\right)|y|^{-\frac{1}{2}}\left(\begin{smallmatrix} |y| & 0 \\ 0 & 1 \end{smallmatrix}\right)\big).
$$
Clearly, this function is well-defined since $|\varepsilon|=1$ for any $\varepsilon\in E_F$.
Let us define the Gamma factor $\Gamma(s,f)$ of $L(s,f)$ as follows:
$$
\Gamma(s,f):=[F_J^\times:F_+^\times]\Big(\frac{\pi}{2}\Big)^{|\Sigma(\mathbb{C})|}\prod_{\tau\in I_F}\frac{\Gamma(s-m_\tau)}{(2\pi)^{s-m_\tau}},
$$
where $F_J^\times:=\{\xi\in F^\times\mid \tau(\xi)>0\ \forall\ \tau\in J\}$ and $F_+^\times:=F_{\Sigma(\mathbb{R})}^\times=F^\times\cap F_{\mathbb{R},+}^\times$.

\begin{proposition}
Let $d^\times y$ be the standard Haar measure on $E_F\backslash F_{\mathbb{R},+}^\times$.
Then, there is an integral representation of the $L$-function as follows:
\begin{equation}\label{integral:representation:lfunction}
D_F^s\Gamma(s,f)L(s,f)=\sum_{a}\int_{E_F\backslash F_{\mathbb{R},+}^\times}f_a(y)|ay|^s d^\times y
\end{equation}
for $\mathrm{Re}(s)>\max \{m_\tau\mid\tau\in I_F\}$,
where $a$ runs through a complete representative set of the narrow class group of $F$.
\end{proposition}
\begin{proof}
Let us denote $E:=E_F\backslash F_{\mathbb{R},+}^\times$ and 
$$
W(y):=\prod_{\tau\in\Sigma(\mathbb{R})}|y_\tau|^{-m_\tau}e^{-2\pi|y_\tau|}\prod_{\tau\in\Sigma(\mathbb{C})}|y_\tau|^{-2m_\tau}K_0(4\pi |y_\tau|),
$$
where $K_j$ is the $j$-th modified Bessel function of the second kind.
From Hida {\cite[Theorem 6.1]{Hida:1994}}, we obtain the following equation when $\mathrm{Re}(s)>\max\{m_\tau\mid\tau\in I_F\} $:
\begin{align*}
&\int_E f_{a}(y)|ay|^sd^\times y=|a|^s\sum_{\alpha\in E_F\backslash F_J^\times}a_f(a\alpha \mathfrak{d}_F)\sum_{\varepsilon\in E_F}\int_E |\varepsilon y|^s W(\alpha \varepsilon y)d^\times \varepsilon y \\
&=\sum_{\alpha}\frac{a_f(a\alpha \mathfrak{d}_F)}{N(aO_F)^s}\int_{ F_{\mathbb{R},+}^\times } |y|^s W(\alpha y)d^\times y=D_F^s\Gamma(s,f)\sum_{\alpha}\frac{a_f( a\alpha\mathfrak{d}_F)}{N( a\alpha\mathfrak{d}_F)^s}.
\end{align*}
Note that  
$$
\sum_{a,\alpha}\frac{a_f( a\alpha\mathfrak{d}_F)}{N( a\alpha\mathfrak{d}_F)^s}=\sum_{\beta\in F_+^\times\backslash F_J^\times}\sum_a\sum_{\alpha\in E_F\backslash F_+^\times}\frac{a_f( a\alpha\beta\mathfrak{d}_F)}{N( a\alpha\beta\mathfrak{d}_F)^s}=\sum_{\beta}\sum_{\mathfrak{a}}\frac{a_f(\mathfrak{a}\beta\mathfrak{d}_F)}{N(\mathfrak{a}\beta\mathfrak{d}_F)^s},
$$
where $\mathfrak{a}$ runs through the ideals of $F$.
Since each $\beta\mathfrak{d}_F$ induces a permutation on the index $\mathfrak{a}$, we obtain the desired result.
\end{proof}

\subsection{Functional equations of {\it L}-functions of cusp forms}\label{subsec:fricke:func:eq}
Let us denote $\varpi_\mathfrak{N}:=\prod_{v|\mathfrak{N}}\varpi_v^{\mathrm{ord}_v(\mathfrak{N})}$.
Let us define a function $W_\mathfrak{N}f:\mathrm{GL}_2(\mathbb{A}_F)\to L_\mathbf{k}$ as follows:
$$
W_{\mathfrak{N}}f(g):=i^{\sum_{\tau\in\Sigma(\mathbb{R})}k_\tau}N(\mathfrak{N})^{1-\frac{k}{2}}\overline{\omega}_\chi(\det g)f\big(g\left(\begin{smallmatrix} 0 & -1 \\ \varpi_\mathfrak{N} & 0 \end{smallmatrix}\right)\big),
$$
where $\omega_\chi:=\chi|\cdot|^{k-2}$ is a unitarization of $\chi$. From the straightforward computation, we observe that $W_{\mathfrak{N}}$ is an involution on the space $S_{(\mathbf{k},\mathbf{m}),J}(\mathfrak{N},\chi)\oplus S_{(\mathbf{k},\mathbf{m}),J}(\mathfrak{N},\overline{\chi})$.
\begin{lemma}\label{functional:equation:lemma}
For a finite idele $a$ of $F$ and an element $y$ of $F_\mathbb{R}^\times$,
$$
N(\mathfrak{N})^{\frac{s}{2}}f_a(y^{-1})|ay^{-1}|^s=W_f^\prime N(\mathfrak{N})^{\frac{k-s}{2}}(W_\mathfrak{N}f)_{\varpi_\mathfrak{N}/a}(y)|(\varpi_\mathfrak{N}/a)y|^{k-s},
$$
where $W_f^\prime:=(-1)^{\sum_{\tau\in\Sigma(\mathbb{R})-J}k_\tau}(-1)^{\sum_{\tau\in\Sigma(\mathbb{C})}(k_\tau-1)}$.
\end{lemma}

\begin{proof}
Let us denote 
$$
g=\left(\begin{smallmatrix} a & 0 \\ 0 & 1 \end{smallmatrix}\right)|y|^{\frac{1}{2}}\left(\begin{smallmatrix} 1/|y| & 0 \\ 0 & 1 \end{smallmatrix}\right),\ g^\prime=\left(\begin{smallmatrix} \varpi_\mathfrak{N}/a & 0 \\ 0 & 1 \end{smallmatrix}\right)|y|^{-\frac{1}{2}}\left(\begin{smallmatrix} |y| & 0 \\ 0 & 1 \end{smallmatrix}\right).
$$
From the straighforward computation, we obtain that 
$$
W_\mathfrak{N}f(g^\prime)=i^{\sum_{\tau\in\Sigma(\mathbb{R})}  k_\tau}N(\mathfrak{N})^{\frac{k}{2}-1}\omega_\chi(a)f\big(\left(\begin{smallmatrix} 0 & -1/a \\ 1 & 0 \end{smallmatrix}\right)|y|^{-\frac{1}{2}}\left(\begin{smallmatrix} |y| & 0 \\ 0 & 1 \end{smallmatrix}\right)\big).
$$
Note that we have the following decomposition:
$$
\left(\begin{smallmatrix} 0 & 1 \\ -1 & 0 \end{smallmatrix}\right)\left(\begin{smallmatrix} 0 & -1/a \\ 1 & 0 \end{smallmatrix}\right)|y|^{-\frac{1}{2}}\left(\begin{smallmatrix} |y| & 0 \\ 0 & 1 \end{smallmatrix}\right)=a^{-1}g\left(\begin{smallmatrix} 0 & 1 \\ -1 & 0 \end{smallmatrix}\right)_\infty,
$$
where $\left(\begin{smallmatrix} 0 & 1 \\ -1 & 0 \end{smallmatrix}\right)\in\mathrm{GL}_2(F)$ and $\left(\begin{smallmatrix} 0 & 1 \\ -1 & 0 \end{smallmatrix}\right)_\infty\in C^+_{F,\infty}$.
In conclude, we obtain that
$$
W_\mathfrak{N}f(g^\prime)=(-1)^{\sum_{\tau\in\Sigma(\mathbb{R})-J}k_\tau}N(\mathfrak{N})^{\frac{k}{2}-1}|a|^{k-2}f(g)\big(\left(\begin{smallmatrix} 0 & 1 \\ -1 & 0 \end{smallmatrix}\right)\mathbf{x}\big).
$$
By using the above equations and the definition of $f_a$, we are done.
\end{proof}

Let us define a completed $L$-function $\Lambda(s,f)$ by 
$$
\Lambda(s,f):=D_F^sN(\mathfrak{N})^{\frac{s}{2}}\Gamma(s,f)L(s,f).
$$
\begin{proposition}\label{Lftn:analytic:conti} 
$\Lambda(s,f)$ has an analytic continuation to entire $s\in\mathbb{C}$, and
\begin{equation}\label{functional:eq}
\Lambda(s,f)=W_f\Lambda(k-s,W_\mathfrak{N} f),
\end{equation}
where 
$
W_f:=(-1)^{|\Sigma(\mathbb{C})|}W_f^\prime=(-1)^{\sum_{\tau\in\Sigma(\mathbb{R})-J}k_\tau}(-1)^{\sum_{\tau\in\Sigma(\mathbb{C})}k_\tau}.
$
\end{proposition}

\begin{proof}
Let us denote by $E:=E_F\backslash F_{\mathbb{R},+}^\times$, $E_1:=\{y\in E:|y|\geq 1\}$, and $E_2:=\{y\in E : |y|\leq 1\}$.
Splitting the integration (\ref{integral:representation:lfunction}) and applying Lemma \ref{functional:equation:lemma} and the change of variable, we can write $\Lambda(s,f)$ as follows:
\begin{align*}
\Lambda(s,f)=&\sum_a N(\mathfrak{N})^{\frac{s}{2}}\Big(\int_{E_1}+\int_{E_2}\Big)f_a(y)|ay|^{s} d^\times y=\sum_a N(\mathfrak{N})^{\frac{s}{2}}\int_{E_1}f_a(y)|ay|^{s} d^\times y \\
&+(-1)^{|\Sigma(\mathbb{C})|}W_f^\prime N(\mathfrak{N})^{\frac{k-s}{2}}\int_{E_1}(W_\mathfrak{N}f)_{\varpi_\mathfrak{N}/a}(y)|(\varpi_\mathfrak{N}/a)y|^{k-s}d^\times y.
\end{align*}
Note that this is well-defined for entire $s\in\mathbb{C}$ due to the behavior of the integrands on $E_1$.
Since $W_\mathfrak{N}^2=1$, $W_f^2=1$, and the map $a\mapsto \varpi_\mathfrak{N}/a$ induces an automorphism on the narrow ray class group of $F$, we are done.
\end{proof}

From this, we can obtain an approximate functional equation: Let us set $|\mathbf{m}|:=\max\{m_\tau\mid \tau\in I_F\}\geq 0$. Let $s\in\mathbb{C}$ such that $\mathrm{Re}(s)>|\mathbf{m}|-2$. Let us denote by
\begin{align*}
V_{j,s}(x)&:=\frac{1}{2\pi i}\int_{2+i\mathbb{R}}\Phi\big((-1)^{j-1}t\big)\Gamma(s+t,f)x^{-t}\frac{dt}{t},
\end{align*} 
where $\Phi(t):=\int_0^\infty\Psi(y)y^t\frac{dy}{y}$, $\Psi$ is a normalized compactly supported $C^\infty$-function. Note that the integrand of $V_{j,s}$ is bounded, decays rapidly as $\mathrm{Im}(t)\rightarrow\infty$, and analytic on $\mathrm{Re}(t)>-1$ except a simple pole at $t=0$. Thus, we obtain the following bounds by shifting the contour of $V_{j,s}
$:
\begin{align}\begin{split}\label{aux.func.esti}
V_{j,s}(x)&=O(\Gamma(s+r,f)x^{-r})\text{ for all }r>0 \text{ as }x\rightarrow\infty, \\
V_{j,s}(x)&=\Gamma(s,f)+O\Big(\Gamma\Big(s-\frac{1}{2},f\Big)x^{\frac{1}{2}}\Big) \text{ as }x\rightarrow 0,
\end{split}\end{align}
where the implicit constants depend only on $\Psi$. 
The residue theorem implies that
\begin{align*}
\Gamma(s,f)L(s,f)=\frac{1}{2\pi i}\int_{2+i\mathbb{R}}-\int_{-\frac{1}{2}+i\mathbb{R}}\Phi(t)\Gamma(s+t,f)L(s+t,f)y^t\frac{dt}{t}.
\end{align*}
Applying the equation (\ref{functional:eq}) to the integrands of the above equation and moving the contour, we obtain  
for $\frac{k}{2}-1<\mathrm{Re}(s)<\frac{k}{2}+1$ that
\begin{align}\label{approx:func:eq:'}
\Gamma(s,f)L(s,f)=\sum_{\mathfrak{a}}\frac{a_f(\mathfrak{a})}{N(\mathfrak{a})^s}V_{1,s}\Big(\frac{N(\mathfrak{a})}{y}\Big) 
+W_f\sum_{\mathfrak{a}}\frac{a_{W_\mathfrak{N} f}(\mathfrak{a})}{N(\mathfrak{a})^{k-s}}V_{2,k-s}\Big(\frac{N(\mathfrak{a})y}{N(\mathfrak{N})}\Big).
\end{align}

\subsection{Gauss sum and twisted cusp form}\label{twisted:cuspforms}

Let $\mathfrak{c}$ be an integral ideal of $F$.
Let $\xi$ be a ray class character modulo $\mathfrak{c}$. For a finite idele $x$ of $F$, let us denote $\xi(x):=\prod_{v\nmid\mathfrak{c}}\xi(v)^{\mathrm{ord}_v(x)}$ by abusing the notation. Then, $\xi$ is a finite order Hecke character modulo $\mathfrak{c}$. 
From now on, by abusing the notation again, let us denote by $\xi(x\mathfrak{a}):=\xi(x)\xi(\mathfrak{a})$ for a finite idele $x$ of $F$ and an ideal $\mathfrak{a}$ of $F$.

Let us denote by $(\mathfrak{c}^{-1}/O_F)^\times:=\{\alpha+O_F\in\mathfrak{c}^{-1}/O_F\mid\mathrm{Ann}_{O_F}(\alpha+O_F)=\mathfrak{c}\}$. 
Let $\mathcal{R}_\mathfrak{c}$ be the image of a complete representative set of $(\mathfrak{c}^{-1}/O_F)^\times$ under the embeddings $F\hookrightarrow \prod_{v|\mathfrak{c}}F_{v}\hookrightarrow \mathbb{A}_F$.
 Note that the map $\mathcal{R}_\mathfrak{c}\rightarrow \prod_{v\mid\mathfrak{c}}(O_v/\mathfrak{c}O_v)^\times$, $u\mapsto \varpi_\mathfrak{c}u$ is a bijection. 
Let us define the {\it Gauss sum} and {\it twisted cusp form}:
\begin{align*}
G(\xi):=\frac{1}{\xi(d_F)}\sum_{u\in\mathcal{R}_\mathfrak{c}}\xi(\varpi_\mathfrak{c} u)\mathbf{e}_F\Big(\frac{u}{d_F}\Big),\ 
f\otimes\xi(g):=\frac{\xi(\det g)}{G(\xi)}\sum_{u\in\mathcal{R}_\mathfrak{c}}\xi(\varpi_\mathfrak{c} u)f\big(g\left(\begin{smallmatrix} 1 & u \\ 0 & 1 \end{smallmatrix}\right)\big),
\end{align*}
where $d_F$ is a finite idele of $F$ corresponding to $\mathfrak{d}_F$, $g\in\mathrm{GL}_2(\mathbb{A}_F)$, and $\xi(\varpi_\mathfrak{c} u):=\xi|_{F^\times_\mathfrak{c}}(\varpi_\mathfrak{c} u)$.
Note that these definitions do not depend on the choice of $d_F$ and $\mathcal{R}_\mathfrak{c}$.
Then, $f\otimes\xi\in S_{(\mathbf{k},\mathbf{m}),J}(\mathfrak{N}\cap\mathfrak{c}^2,\chi\xi^2)$ and $a_{f\otimes\xi}(\mathfrak{a})=a_f(\mathfrak{a})\xi(\mathfrak{a})$ by Hida \cite[Section 6]{Hida:1994}.
 We have the following relation between the involution $W_\mathfrak{N}$ and twisting $\otimes\xi$:

\begin{proposition}\label{fricke:twist:commute} If $\xi$ is primitive and $\mathfrak{c}$ is coprime to $\mathfrak{N}$, then
$$
W_{\mathfrak{N}\mathfrak{c}^2}(f\otimes\xi)=\omega_\chi(\varpi_\mathfrak{c})W(\xi)(W_{\mathfrak{N}}f)\otimes\overline{\xi},
$$
where $W(\xi):=N(\mathfrak{c})^{-1} G(\overline{\xi})^2 \xi(\varpi_{\mathfrak{c}^2\mathfrak{N}})$, a complex number of the absolute value $1$.
\end{proposition}

\begin{proof} 
It is clear from the definition that
\begin{align*}
W_{\mathfrak{N}\mathfrak{c}^2}(f\otimes\xi)(g)&=\frac{N(\mathfrak{c})^{2-k}\xi(\varpi_{\mathfrak{N}\mathfrak{c}^2})C}{\omega_\chi\xi(\det g)G(\xi)}\sum_{u\in\mathcal{R}_\mathfrak{c}}\xi(\varpi_\mathfrak{c} u)f\Big(g\big(\begin{smallmatrix} 0 & -1 \\ \varpi_{\mathfrak{N}\mathfrak{c}^2} & 0 \end{smallmatrix}\big)\left(\begin{smallmatrix} 1 & u \\ 0 & 1 \end{smallmatrix}\right)\Big), \\
(W_{\mathfrak{N}}f)\otimes\overline{\xi}(g)&=\frac{C}{\omega_\chi\xi(\det g)G(\overline{\xi})} \sum_{u\in\mathcal{R}_\mathfrak{c}}\overline{\xi}(\varpi_\mathfrak{c} u)f\Big(g\left(\begin{smallmatrix} 1 & u \\ 0 & 1 \end{smallmatrix}\right)\big(\begin{smallmatrix} 0 & -1 \\ \varpi_\mathfrak{N} & 0 \end{smallmatrix}\big)\Big),
\end{align*}
where $C:=i^{\sum_{\tau\in\Sigma(\mathbb{R})}k_\tau}N(\mathfrak{N})^{1-\frac{k}{2}}$.
For each $u\in\mathcal{R}_\mathfrak{c}$, there exists $u^\prime\in \mathcal{R}_\mathfrak{c}$ and $a_u\in O_\mathfrak{c}:=\prod_{v\mid\mathfrak{c}}O_v$ such that $\varpi_\mathfrak{c} u \varpi_\mathfrak{c} u^\prime=\varpi_\mathfrak{c}a_u-1$ in $\prod_{v\mid\mathfrak{c}}O_v$.
Let $A_u$ be an element of $\mathrm{GL}_2(\mathbb{A}_F)$ such that
$(A_u)_\mathfrak{c}=\left(\begin{smallmatrix} -\varpi_\mathfrak{c}u^\prime & -a_u \\ \varpi_\mathfrak{c} & \varpi_\mathfrak{c}u \end{smallmatrix}\right)\in\mathrm{SL}_2(O_\mathfrak{c})$ and $(A_u)_v=\left(\begin{smallmatrix} 1 & 0 \\ 0 & 1\end{smallmatrix}\right)$ if $v\nmid\mathfrak{c}$ . 
Then, by comparing the both sides at each place, we observe that
$$
\big(\begin{smallmatrix} 0 & -1 \\ \varpi_{\mathfrak{N}\mathfrak{c}^2} & 0 \end{smallmatrix}\big)\left(\begin{smallmatrix} 1 & u \\ 0 & 1 \end{smallmatrix}\right)=\varpi_\mathfrak{c}\left(\begin{smallmatrix} 1 & u^\prime \\ 0 & 1 \end{smallmatrix}\right)\left(\begin{smallmatrix} 0 & -1 \\ \varpi_\mathfrak{N} & 0 \end{smallmatrix}\right)A_u .
$$
From this equation, we can rewrite $W_{\mathfrak{Nc^2}}(f\otimes\xi)$ as follows:
\begin{align*}
W_{\mathfrak{N}\mathfrak{c}^2}(f\otimes\xi)(g)=\frac{\xi(\varpi_{\mathfrak{N}\mathfrak{c}^2})\xi(-1_\infty)\omega_\chi(\varpi_\mathfrak{c})C}{\omega_\chi\xi(\det(g))G(\xi)}\sum_{u\in\mathcal{R}_\mathfrak{c}}\overline{\xi}(\varpi_\mathfrak{c} u^\prime)f\Big(g\left(\begin{smallmatrix} 1 & u^\prime \\ 0 & 1 \end{smallmatrix}\right)\big(\begin{smallmatrix} 0 & -1 \\ \varpi_{\mathfrak{N}} & 0 \end{smallmatrix}\big)\Big).
\end{align*} 
Note that the map $u\mapsto u^\prime$ sends $\mathcal{R}_\mathfrak{c}$ to another choice of $\mathcal{R}_\mathfrak{c}$.
Using the properties of Gauss sum (Neukirch, \cite[Chapter VII, Proposition 7.5]{neukirch2013algebraic}) and comparing the above equations, we obtain the desired result.
\end{proof}


\section{Averages of ray class characters and {\it L}-values}\label{subsec:hecke:char}
In this section, we study averages of Hecke characters and $L$-values, which play a crucial role to obtain the main result.

\subsection{Averages of ray class characters}

Let $\mathfrak{p}$ be a prime ideal of $F$ lying above an odd rational prime $p$.
Let us denote by $\mathrm{Cl}_F(\mathfrak{p}^m)$ the ray class group of $F$ of modulus $\mathfrak{p}^m $.
Let us set
$$
\mathrm{Cl}_F(\mathfrak{p}^\infty):=\varprojlim_{m>0} \mathrm{Cl}_F(\mathfrak{p}^m).
$$
Let us denote by $\mathrm{Cl}_F$ the ideal class group of $F$.
From the fundamental exact sequence of ray class groups, we observe that the following sequence is exact:
\begin{equation}\label{fund:exact:seq}
\begin{tikzcd} 
O_\mathfrak{p}^\times \arrow{r} & \mathrm{Cl}_F(\mathfrak{p}^\infty) \arrow{r} & \mathrm{Cl}_F \arrow{r} & 1. 
\end{tikzcd}
\end{equation}

\begin{definition}\label{coates:wiles:extension}
Let us denote by $\Xi_\mathfrak{p}$ the set of $p$-power order primitive ray class characters of $\mathfrak{p}$-power modulus.
\end{definition}
Let $\phi$ be an element of $\Xi_\mathfrak{p}$.
Note that the sign of $\phi$ is totally even so that $\phi$ can be considered as a character of the pro $p$-part of $\mathrm{Cl}_F(\mathfrak{p}^\infty)$.
From now on, we assume that the residue degree and the ramification degree of $\mathfrak{p}$ over $p$ are $1$, so that 
$\mathrm{rank}_{\mathbb{Z}_p}\big(\mathrm{Cl}_F(\mathfrak{p}^\infty)\big)\leq 1$.
We only consider the case that $\mathrm{rank}_{\mathbb{Z}_p}\big(\mathrm{Cl}_F(\mathfrak{p}^\infty)\big)=1$. If the rank is zero, then $\Xi_\mathfrak{p}$ is a finite set so that our main theorem vacuously true.


Let $n_0>0$ be an integer. Let $n$ be the exponent of the conductor of $\phi$. 
Let us assume that $n>n_0$. 
As $b\mapsto\phi(1+\varpi_\mathfrak{p}^{n-n_0}b)$ is an additive character of $O_\mathfrak{p}/\mathfrak{p}^{n_0}O_\mathfrak{p}$, there is an element $c\in O_\mathfrak{p}$ such that 
$
\phi(1+\varpi_\mathfrak{p}^{n-n_0}b)=\mathbf{e}_F(bc/\varpi_\mathfrak{p}^{n_0})
$
for any $b\in O_\mathfrak{p}$. 
Let us denote $c_\phi:=c$.
Note that $c_\phi\in O_\mathfrak{p}^\times$ as $\phi$ is primitive.
Let us define a {\it partial Gauss sum} as follows:
$$
G_1(\phi):=\sum_{a\equiv c_\phi(\mathfrak{p}^{n_0})}\phi(a)\mathbf{e}_F\Big(\frac{a}{\varpi_\mathfrak{p}^n }\Big),
$$
where $a$ runs through $(O_\mathfrak{p}/\mathfrak{p}^nO_\mathfrak{p})^\times$. Note that this definition does not depend on the choice of $a$'s. Note that $G(\phi)=\sum_{a}\phi(a)\mathbf{e}_F\big(a/\varpi_\mathfrak{p}^n \big)$ as $\mathfrak{p}$ is unramified in $F$.
\begin{definition}\label{galois:averages}
For an ideal $\mathfrak{a}$ of $F$, let us define the averages as follows: 
\begin{align*}
G_{\mathrm{av}}(\phi,\mathfrak{a})&:=\frac{1}{|G|p^n}\sum_{\sigma\in G}G_1(\overline{\phi}^\sigma)G(\phi^\sigma)\phi^\sigma(\mathfrak{a}),
\\
G_{\mathrm{av}}^\iota(\phi,\mathfrak{a})&:=\frac{1}{|G|p^n}\sum_{\sigma\in G}G_1(\overline{\phi}^\sigma)G(\phi^\sigma)W(\phi^\sigma)\overline{\phi}^\sigma(\mathfrak{a}),
\end{align*}
where $W(\phi):=p^{-n} G(\overline{\phi})^2 \phi(\varpi_{\mathfrak{p}^{2n}\mathfrak{N}})$ and $G=\mathrm{Gal}\big(\mathbb{Q}(\phi)/\mathbb{Q}(\mu_{p^{n_0}})\big)$.
Let us define 
$$
\phi_{\mathrm{av}}(\mathfrak{a}):=\frac{1}{|G|}\sum_{\sigma\in G}\phi^\sigma(\mathfrak{a}).
$$
Note that the averages vanish if $\mathfrak{p}|\mathfrak{a}$.
\end{definition}

\begin{lemma}\label{gal:av:nonzero} 
For an ideal $\mathfrak{a}$ of $F$, $G_{\mathrm{av}}(\phi,\mathfrak{a})=\phi_{\mathrm{av}}(\mathfrak{a})$.
\end{lemma}

\begin{proof}  
From the definition, we observe that
\begin{align*}
G_{\mathrm{av}}(\phi,\mathfrak{a})=&\frac{1}{p^n}\sum_{a_1,a}\phi_{\rm av}(a_1^{-1}a\mathfrak{a})\mathbf{e}_F\Big(\frac{a_1+a}{\varpi_\mathfrak{p}^n}\Big)=\frac{1}{p^n}\sum_{a_1,a_2}\phi_{\rm av}(a_2\mathfrak{a})\mathbf{e}_F\Big(\frac{a_1(1+a_2)}{\varpi_\mathfrak{p}^n}\Big)\\
=&\frac{1}{p^{n_0}}\sum_{a_2\equiv -1(\mathfrak{p}^{n-n_0})}\phi_{\rm av}(a_2\mathfrak{a})\mathbf{e}_F\Big(\frac{c_\phi(1+a_2)}{\varpi_\mathfrak{p}^n}\Big),
\end{align*}
where $a_1\equiv c_\phi(p^{n_0})$, $a$ and $a_2$ runs through $(O_\mathfrak{p}/\mathfrak{p}^nO_\mathfrak{p})^\times$. The last equality holds 
due to the parametrization $a_1=c_\phi+\varpi_\mathfrak{p}^{n_0}b$, $b\in O_\mathfrak{p}/\mathfrak{p}^{n-n_0}O_\mathfrak{p}$ and the orthoginality of characters.
Note that $\phi_{\rm av}(a_2\mathfrak{a})=\phi(a_2)\phi_{\rm av}(\mathfrak{a})$ and $\phi(a_2)=\mathbf{e}_F(-c_\phi(1+a_2)/\varpi_\mathfrak{p}^n)$, so we are done.
\end{proof}

We may take a complete representative set $\{\mathfrak{a}_i^{-1}\}_{i=1}^{h_F}$ of $\mathrm{Cl}_F$ such that $\mathfrak{a}_1=O_F$ and each $\mathfrak{a}_i$ is an integral ideal coprime to $\mathfrak{p}$.
If $\mathfrak{a}$ is an ideal of $F$ coprime to $\mathfrak{p}$, then $\mathfrak{a}=\alpha\mathfrak{a}_j^{-1}$ for some $\alpha\in F^\times\cap O_\mathfrak{p}^\times$ and some $j$ due to the exact sequence (\ref{fund:exact:seq}).
Let $W_\mathfrak{p}$ be the split image of $(O_F/\mathfrak{p})^\times$ in $O_\mathfrak{p}^\times$. Note that $O_\mathfrak{p}^\times\cong W_\mathfrak{p}\times (1+\mathfrak{p}O_\mathfrak{p})$ and $|W_\mathfrak{p}|=p-1$.
Let us assume that $p$ is coprime to $h_F$, so that the map $x\mapsto x^{h_F(p-1)}:O_\mathfrak{p}^\times\rightarrow 1+\mathfrak{p}O_\mathfrak{p}$ is $p-1$ to 1 surjective. 
Let $\alpha_j$ be an element of $F^\times\cap(1+\mathfrak{p}O_\mathfrak{p})$ such that $\mathfrak{a}_j^{h_F(p-1)}=\alpha_j O_F$.
Let $b_j$ be the element of $1+\mathfrak{p}O_\mathfrak{p}$ such that $b_j^{h_F(p-1)}=\alpha_j$. 
Then, the inverse image of $\alpha_j$ under the map $x\mapsto x^{h_F(p-1)}$ is given by $\{\kappa b_j\}_{\kappa\in W_\mathfrak{p}}$.
We have the following useful lemmas:

\begin{lemma}\label{useful:lemma}
Let $\mathfrak{a}$ be an ideal of $F$ such that $\phi_{\mathrm{av}}(\mathfrak{a})\neq 0$, so that $\mathfrak{a}=\alpha\mathfrak{a}_j^{-1}$ for some $\alpha\in F^\times\cap O_\mathfrak{p}^\times$ and $j$. Then, $\phi_{\mathrm{av}}(\mathfrak{a})=\phi(\mathfrak{a})$ and $\alpha^{h_F(p-1)}\in \alpha_j+\mathfrak{p}^{n-n_0}O_\mathfrak{p}$, or equivalently,
$$
\alpha\in\bigcup_{\kappa\in W_\mathfrak{p}} \kappa b_j+\mathfrak{p}^{n-n_0}O_\mathfrak{p}.
$$
The converse also holds.
\end{lemma}

\begin{proof}
The restriction of $\phi$ gives the following isomorphism of cyclic groups:
$$
\phi:\frac{1+\mathfrak{p}O_\mathfrak{p}}{1+\mathfrak{p}^nO_\mathfrak{p}}\cong\mu_{p^{n-1}}.
$$
Thus, it is clear that $\phi_{\mathrm{av}}(\mathfrak{a})\neq 0$ if and only if $\phi(\alpha^{h_F(p-1)}\alpha_j^{-1})=\phi(\mathfrak{a})^{h_F(p-1)}\in\mu_{p^{n_0}}$ as $p\nmid h_F(p-1)$.
If $\phi(\mathfrak{a})^{h_F(p-1)}\in\mu_{p^{n_0}}$, then $\alpha^{h_F(p-1)}\in\alpha_j+\mathfrak{p}^{n-n_0}O_\mathfrak{p}$.
The converse direction is clear. So we are done.
\end{proof}

\begin{lemma}\label{useful:lemma:2}
For elements $c,d\in O_\mathfrak{p}^\times$ and an integer $n$ with $\lceil\frac{n}{2}\rceil\geq n_0$, 
$$
\frac{1}{p^{n}}\sum_{a\equiv c(\mathfrak{p}^{n_0})}\mathbf{e}_F\Big(\frac{a+a^{-1}d}{\varpi_\mathfrak{p}^n}\Big)\leq p^{2-\frac{n}{2}}.
$$
\end{lemma}
\begin{proof}
Let us write $m=\lceil\frac{n}{2}\rceil$ and $P(t)=t+t^{-1}d$. By abusing the notation, let us decompose $a=ac$, where $a\in\frac{1+\mathfrak{p}^{n_0}O_\mathfrak{p}}{1+\mathfrak{p}^nO_\mathfrak{p}}$.  Then, we can write
$$
\sum_{a\equiv c(\mathfrak{p}^{n_0})}\mathbf{e}_F\Big(\frac{a+a^{-1}d}{\varpi_\mathfrak{p}^n}\Big)=\sum_{a\equiv 1(\mathfrak{p}^{n_0})}\mathbf{e}_F\Big(\frac{P(ac)}{\varpi_\mathfrak{p}^n}\Big)=\sum_{b,r}\mathbf{e}_F\Big(\frac{P(a_{b,r}c)}{\varpi_\mathfrak{p}^n}\Big),
$$
where $a_{b,r}=r+\varpi_\mathfrak{p}^{m}b$, $r\in\frac{1+\mathfrak{p}^{n_0}O_\mathfrak{p}}{1+\mathfrak{p}^mO_\mathfrak{p}}$, $b\in O_\mathfrak{p}/\mathfrak{p}^{n-m}O_\mathfrak{p}$. Note that $P(a_{b,r}c)\equiv P(cr)+\varpi_\mathfrak{p}^{m}b(c-c^{-1}dr^{-2})\ (\mathfrak{p}^n)$. Therefore, by using the orthogonality of characters, we rewrite the above equation as follows:
$$
\sum_{r}\mathbf{e}_F\Big(\frac{P(cr)}{\varpi_\mathfrak{p}^n}\Big)\sum_{b(\mathfrak{p}^{n-m})}\mathbf{e}_F\Big(\frac{b(c-r^{-2}c^{-1}d)}{\varpi_\mathfrak{p}^{n-m}}\Big)=p^{n-m}\sum_{r^2\equiv c^{-2}d(\mathfrak{p}^{n-m})}\mathbf{e}_F\Big(\frac{P(cr)}{\varpi_\mathfrak{p}^n}\Big).
$$
Note that $r\mapsto r^2$ is a group automorphism on $\frac{1+\mathfrak{p}^{n_0}O_\mathfrak{p}}{1+\mathfrak{p}^mO_\mathfrak{p}}$. Also note that $n-m=m$ or $m-1$. So we are done.
\end{proof}

\begin{proposition}\label{root:num:gal:av} Let $\mathfrak{a}$ be an ideal of $F$. If $\mathfrak{p}\nmid\mathfrak{N}$ and $n>2n_0$, then
$$
|G_{\mathrm{av}}^\iota(\phi,\mathfrak{a})|\leq (p-1)p^{n_0+2-\frac{n}{2}}.
$$
\end{proposition}

\begin{proof}
We may assume that $\mathfrak{p}\nmid\mathfrak{a}$. Then, $\mathfrak{a}=\alpha\mathfrak{a}_j^{-1}$ for some $\alpha\in O_F\cap O_\mathfrak{p}^\times$ and $j$. 
Let $c_n\in\mathbb{Z}$ be a mod $p^n$ inverse of $h_F$. 
Note that $\mathfrak{N}^{h_F}=\beta O_F$ for some $\beta\in O_F\cap O_\mathfrak{p}^\times$. Then, ${\phi}(\beta^{c_n} O_F)=\phi(\mathfrak{N})^{h_Fc_n}=\phi(\mathfrak{N})=\phi(\varpi_{\mathfrak{p}^{2n}\mathfrak{N}})$ as $\phi(\varpi_\mathfrak{p})=\phi(O_F)$.
Hence, we obtain the following from the definition and Lemma \ref{useful:lemma}:
\begin{align*}
G_{\mathrm{av}}^\iota(\phi,\mathfrak{a})=&\frac{1}{p^n}\sum_{a_1,a}\overline{\phi}_{\rm av}(a_1a c\mathfrak{a}_j^{-1})\mathbf{e}_F\Big(\frac{a_1+a}{\varpi_\mathfrak{p}^n}\Big)=\frac{1}{p^n}\sum_{a_1,a_2}\overline{\phi}_{\mathrm{av}}(a_2c\mathfrak{a}^{-1}_j)\mathbf{e}_F\Big(\frac{a_1+a_1^{-1}a_2}{\varpi_\mathfrak{p}^n}\Big) \\
=&\sum_{\kappa\in W_\mathfrak{p}}\sum_{a_2\in\kappa b_j c^{-1}+\mathfrak{p}^{n-n_0}O_\mathfrak{p}}\overline{\phi}(a_2c\mathfrak{a}_j^{-1})\frac{1}{p^n}\sum_{a_1}\mathbf{e}_F\Big(\frac{a_1+a_1^{-1}a_2}{\varpi_\mathfrak{p}^n}\Big),
\end{align*}
where $a_1\equiv c_\phi(p^{n_0})$, $a$ and $a_2$ runs through $(O_\mathfrak{p}/\mathfrak{p}^nO_\mathfrak{p})^\times$, and $c:=\alpha\beta^{-c_n}\in F\cap O_\mathfrak{p}^\times$. 
Appying Lemma \ref{useful:lemma:2} to the above equation, we obtain the desired bound.
\end{proof}

\subsection{Coherent cone decomposition and some estimations}
For a $\mathbb{Q}$-basis $V$ of $F$, let us set 
$$
\mathbb{Z}_+\langle V \rangle:=\Big\{\sum_{a\in V}n_a a \mid n_a\in \mathbb{Z}_+\Big\}\subset F_\mathbb{R}^\times,
$$ 
which is called $\mathbb{Z}$-{\it cone} generated by $V$. 
A $\mathbb{Z}$-cone $\mathbb{Z}_+\langle V \rangle$ is said to be {\it coherent} if we can choose $\rho_\tau\in\{\operatorname{Re}\tau,\operatorname{Im}\tau\}$ for each $\tau\in I_F$ such that the numbers $\rho_\tau(a)$, $a\in V$ are non-zero and of the same sign.

Let $C_F^0$ be a fundamental domain of $F_\mathbb{R}^\times/O_F^\times$, where $O_F^\times$ acts on $F_\mathbb{R}^\times$ via the natural embedding $F^\times\hookrightarrow F_\mathbb{R}^\times$.
Let us denote by $\mathcal{B}$ a finite set of coherent $\mathbb{Z}$-cones which covers $O_F\cap C_F^0$. The existence of such $\mathcal{B}$ is guaranteed by Rohrlich \cite[Proposition 4]{Rohrlich:1989}.
Then, we obtain the following estimations:

\begin{lemma}\label{lattice:esti:2} 
Let $\mathfrak{c}$ be an ideal of $F$ and $\gamma$ an element of $F$. If $\alpha$ is an element of $O_F\cap C_F^0$ such that $\alpha^m-\gamma\in\mathfrak{c}-\{0\}$ for some integer $m>0$, then 
$$
|N(\alpha)|\gg_{F,\gamma,m} N(\mathfrak{c})^{\frac{1}{m}}.
$$
\end{lemma}

\begin{proof} 
Note that $\alpha\in B$ for some $B=\mathbb{Z}_+\langle V\rangle\in\mathcal{B}$. 
We may write $\alpha=\sum_{a\in V}n_a a$. For an element $\tau$ of $I_F$, 
$$
|\alpha^\tau|\geq |\rho_\tau(\alpha)|=\sum_{a\in V} n_a|\rho_\tau(a)|\gg_{B,\tau} 1
$$ 
as $B$ is coherent. In conclude, we obtain the following bound:
$$
N(\mathfrak{c})\leq |N(\alpha^m-\gamma)|\leq |N(\alpha)|^m\prod_{\tau\in I_F}\Big( 1+\frac{|\gamma^\tau|}{|\alpha^\tau|^{m}}\Big)\ll_{B,\gamma,I_F,m} |N(\alpha)|^m.
$$
Note that $\mathcal{B}$ depends only on $F$. So we are done.
\end{proof}

\begin{proposition}\label{1st:gal:av:Lvalue:error:1:lemma}
Let $\mathfrak{a},\mathfrak{r}$ be distinct integral ideals of $F$ coprime to $\mathfrak{p}$ so that $\mathfrak{a}=\alpha\mathfrak{a}_j^{-1}$, $\mathfrak{r}=\gamma\mathfrak{a}_l^{-1}$ for some elements $\alpha,\gamma\in O_F\cap O_\mathfrak{p}^\times$ and some $j,l$.
If $G_{\mathrm{av}}(\phi,\mathfrak{r}^{-1}\mathfrak{a})\neq 0$, then
$$
\alpha\in\bigcup_{\kappa\in W_\mathfrak{p}} (\kappa\beta_{n}+\mathfrak{p}^{n-n_0})\cap C_F^0,\text{ and } |N(\alpha)|>c_{\mathfrak{r}} p^{\frac{n-n_0}{h_F(p-1)}}
$$
for some constant $c_{\mathfrak{r}}>0$ independent on $n$, where
$\beta_{n}\in O_F\cap O_\mathfrak{p}^\times$ is a representative of $b_l^{-1}b_j$ mod $\mathfrak{p}^{n-n_0}$.
\end{proposition}

\begin{proof} 
We may assume that $\alpha\in C_F^0$ since $\mathfrak{a}=\alpha\mathfrak{a}_j^{-1}$ is integral.
Note that $G_{\mathrm{av}}(\phi,\mathfrak{r}^{-1}\mathfrak{a})=\phi_{\mathrm{av}}(\mathfrak{r}^{-1}\mathfrak{a})$ by Lemma \ref{gal:av:nonzero}.
As $\alpha$ and $\alpha_l$ are integral, we observe from Lemma \ref{useful:lemma} that $\alpha^d\alpha_l\in\gamma^d\alpha_j+\mathfrak{p}^{n-n_0}$, where $d=h_F(p-1)$.
If $\alpha^d\alpha_l=\gamma^{d}\alpha_j$, then $(\mathfrak{r}^{-1}\mathfrak{a})^d=O_F$, which implies that $\mathfrak{a}=\mathfrak{r}$. Thus, $\alpha^d-\gamma^d\alpha_l^{-1}\alpha_j\in\alpha_l^{-1}\mathfrak{p}^{n-n_0}-\{0\}$ so that
$$
|N(\alpha)|\gg_{d,F,\gamma,N(\mathfrak{a}_l),p} p^{\frac{n-n_0}{d}} 
$$
by Lemma \ref{lattice:esti:2}.
So we are done.
\end{proof}

\subsection{Average of {\it L}-values}
Let $f\in S_{(\mathbf{k},\mathbf{m}),J}(\mathfrak{N},\chi)$ be a Hecke eigenform.
\begin{definition}\label{galois:average:Lvalue}
For an integral ideal $\mathfrak{r}$ of $F$ coprime to $\mathfrak{p}$, let us define
$$
L_{\mathrm{av}}(f,\phi,\mathfrak{r}):=\frac{1}{|G|p^n}\sum_{a\equiv c_\phi(\mathfrak{p}^{n_0})}\mathbf{e}_F\Big(\frac{a}{\varpi_\mathfrak{p}^n}\Big)\sum_{\sigma\in G}\overline{\phi}^\sigma(a\mathfrak{r})G(\phi^\sigma)L\Big(\frac{k}{2},f\otimes\phi^\sigma\Big),
$$
where $G=\mathrm{Gal}\big(\mathbb{Q}(\phi)/\mathbb{Q}(\mu_{p^{n_0}})\big)$, and $a$ runs through $(O_\mathfrak{p}/\mathfrak{p}^nO_\mathfrak{p})^\times$.
\end{definition} 
Let us assume that $\mathfrak{p}$ is coprime to $\mathfrak{N}$.
Then, by the equation (\ref{approx:func:eq:'}) and Proposition \ref{fricke:twist:commute}, we can write
$L_{\mathrm{av}}(f,\phi,\mathfrak{r})=L_{\mathrm{av},1}(f,\phi,\mathfrak{r})+L_{\mathrm{av},2}(f,\phi,\mathfrak{r})$, where
\begin{align*}
L_{\mathrm{av},1}(f,\phi,\mathfrak{r}):=&\frac{1}{\Gamma(\frac{k}{2},f)}\sum_{\mathfrak{a}}\frac{a_f(\mathfrak{a})G_{\mathrm{av}}(\phi,\mathfrak{r}^{-1}\mathfrak{a})}{N(\mathfrak{a})^{\frac{k}{2}}}V_{1,\frac{k}{2}}\Big(\frac{N(\mathfrak{a})}{y}\Big),  \\
L_{\mathrm{av},2}(f,\phi,\mathfrak{r}):=&\frac{\omega_\chi(\varpi_\mathfrak{p}^n)W_f}{\Gamma(\frac{k}{2},f)}\sum_{\mathfrak{a}}\frac{a_{W_\mathfrak{N} f}(\mathfrak{a})G^\iota_{\mathrm{av}}(\phi,\mathfrak{r}\mathfrak{a})}{N(\mathfrak{a})^{\frac{k}{2}}}V_{2,\frac{k}{2}}\Big(\frac{N(\mathfrak{a})y}{N(\mathfrak{N})p^{2n}}\Big),
\end{align*} 
where $\mathfrak{a}$ runs through the integral ideals.
By Blomer-Brumley \cite{blomer2011ramanujan} and Nakasuji \cite{nakasuji2012generalized},
\begin{equation}\label{rama:peter:bdd}
|a_f(\mathfrak{a})|\leq 2d(\mathfrak{a})N(\mathfrak{a})^{\frac{k-1}{2}+\frac{7}{64}}\ll_{F,\varepsilon} N(\mathfrak{a})^{\frac{k}{2}-\frac{25}{64}+\varepsilon}
\end{equation}
for $\varepsilon>0$ and each integral ideal $\mathfrak{a}$ of $F$,
where $d(\mathfrak{a})$ is the number of the integral ideals of $F$ dividing $\mathfrak{a}$.
Note that the last inequality of the equation (\ref{rama:peter:bdd}) holds as $d(m)\ll_{\varepsilon,F} m^\varepsilon$, which comes from the convergence of the Dedekind zeta function of $F$.
Let us estimate $L_{\mathrm{av},i}(f,\phi,\mathfrak{r})$:

\begin{proposition}\label{1st:gal:av:Lvalue:error:1} 
Let $\mathfrak{r}$ be an integral ideal of $F$ coprime to $\mathfrak{p}$.
Let $c_\mathfrak{r}$ be the constant defined in Proposition \ref{1st:gal:av:Lvalue:error:1:lemma}, which does not depend on $n$.
Let $y>p^n$ be a real number. For a small $\varepsilon>0$ and all sufficiently large $n$,
\begin{align*}
L_{\mathrm{av},1}(f,\phi,\mathfrak{r})-a_f(\mathfrak{r})
\ll_{c_\mathfrak{r},\varepsilon,F,k,\{N(\mathfrak{a}_i)\}_i,p} \frac{1}{y^{\frac{1}{2}}}+\frac{1}{p^{\frac{n-n_0}{h_F(p-1)}(\frac{25}{64}-\varepsilon)}}+\frac{y^{\frac{39}{64}+\varepsilon}}{p^n}.
\end{align*}
\end{proposition}

\begin{proof}
From the estimations (\ref{aux.func.esti}), (\ref{rama:peter:bdd}) and Proposition \ref{1st:gal:av:Lvalue:error:1:lemma}, we observe that
\begin{align*}
&L_{\mathrm{av},1}(f,\phi,\mathfrak{r})-a_f(\mathfrak{r})=\frac{a_f(\mathfrak{r})}{\Gamma(\frac{k}{2},f)}V_{1,\frac{k}{2}}\Big(\frac{1}{y}\Big)-a_f(\mathfrak{r})+\sum_{\mathfrak{a}\neq\mathfrak{r}}\ll_{\varepsilon,k,\{N(\mathfrak{a}_i)\}_i} \frac{1}{y^{\frac{1}{2}}}+(^*)+(^{**}), \\
&(^*)=\sum_{i=1}^{h_F}\sum_{\kappa\in W_\mathfrak{p} }\sum_{c_\mathfrak{r} p^{\frac{n-n_0}{d}}<|N(\alpha)|\leq y}\frac{y}{|N(\alpha)|^{\frac{25}{64}-\varepsilon}},\
(^{**})=\sum_{i=1}^{h_F}\sum_{\kappa\in W_\mathfrak{p} }\sum_{|N(\alpha)|>y}\frac{y}{|N(\alpha)|^{\frac{25}{64}-\varepsilon}},
\end{align*}
where $\alpha$ runs through $(\kappa\beta_{n}+\mathfrak{p}^{n-n_0})\cap C_F^0$, $\beta_n$ is a number defined in Proposition \ref{1st:gal:av:Lvalue:error:1:lemma}, and $d=h_F(p-1)$.
For $a\in O_\mathfrak{p}^\times$ and $x>0$, let us set 
$$
u_{a,n}(x):=\#\{\alpha\in (a+\mathfrak{p}^{n-n_0})\cap C_F^0:|N(\alpha)|=x\}.
$$
Then, we obtain from the proof of Rohrlich \cite[Proposition 5]{Rohrlich:1989} that 
$$
\sum_{m\leq x} u_{a,n}(m)\ll \mathrm{max}\Big\{\frac{x}{p^n},1\Big\},
$$ 
where the implicit constant does not depend on $a$ and $n$. Thus, by the Abel summation formula, we obtain that
\begin{align*}
&(^*)\ll_{d,F,p} \sum_{c_\mathfrak{r} p^{\frac{n-n_0}{d}}<m\leq y} \frac{u_{\kappa\beta_{n},n}(m) }{m^{\frac{25}{64}-\varepsilon}}\ll\frac{1}{p^{\frac{n-n_0}{d}(\frac{25}{64}-\varepsilon)}}+\frac{y^{\frac{39}{64}+\varepsilon}}{p^n} \\
&+\int_{c_\mathfrak{r} p^{\frac{n-n_0}{d}}}^{p^n}\frac{dx}{x^{\frac{89}{64}-\varepsilon }}+\int_{p^n}^{y}\frac{dx}{p^nx^{\frac{25}{64}-\varepsilon }}dx\ll_{c_\mathfrak{r},\varepsilon,F} \frac{1}{p^{\frac{n-n_0}{d}(\frac{25}{64}-\varepsilon)}}+\frac{y^{\frac{39}{64}+\varepsilon}}{p^n},
\end{align*}
\begin{align*}
(^{**})\ll_{d,F,p}  \sum_{m>y}\frac{u_{\kappa\beta_{n},n}(m)y}{m^{\frac{89}{64}-\varepsilon}}
\ll \frac{y^{\frac{39}{64}+\varepsilon}}{p^n}+\int_{y}^{\infty}\frac{ydx}{p^n x^{\frac{89}{64}-\varepsilon}}\ll_{\varepsilon,F} \frac{y^{\frac{39}{64}+\varepsilon}}{p^{n}}.
\end{align*}
So we are done.
\end{proof}

\begin{proposition}\label{1st:gal:av:Lvalue:error:2} Let $0<y<p^{2n}$ be a real number. For a small $\varepsilon>0$ and all sufficiently large $n$,
\begin{align*}
L_{\mathrm{av},2}(f,\phi,\mathfrak{r})\ll_{\varepsilon,n_0,N(\mathfrak{N}),p} & \frac{p^{(\frac{23}{32}+\varepsilon)n}}{y^{\frac{39}{64}+\varepsilon}}.
\end{align*}
\end{proposition}

\begin{proof}
Note that $W_\mathfrak{N} f$ is also a Hecke eigenform (see Miyake \cite[(1.8)]{miyake1971onautomorphic}). Due to the bound $\sum_{N(\mathfrak{a})=m}1\leq d(m)\ll_\varepsilon m^\varepsilon$, the equation (\ref{rama:peter:bdd}), and Proposition \ref{root:num:gal:av}, 
\begin{align*}\label{1st:gal:av:Lvalue:error:2:esti}
L_{\mathrm{av},2}(f,\phi,\mathfrak{r})&\ll_{\varepsilon,n_0,p}  \sum_{m>0} \frac{m^{-\frac{25}{64}+\varepsilon}}{p^{\frac{n}{2}}}V_{2,\frac{k}{2}}\Big(\frac{my}{N(\mathfrak{N})p^{2n}}\Big)=I+II,\text{ where }\\
&I=\sum_{m>N(\mathfrak{N})p^{2n}/y},\ II=\sum_{0<m<N(\mathfrak{N})p^{2n}/y}.
\end{align*}
Applying the equation (\ref{aux.func.esti}) to the above equation, we obtain that
\begin{align*}
I&\ll_{N(\mathfrak{N})} \sum_{m>N(\mathfrak{N})p^{2n}/y}\frac{p^{\frac{3n}{2}}}{m^{\frac{89}{64}-\varepsilon}y} \leq \frac{p^{\frac{3n}{2}}}{y}\int_{ N(\mathfrak{N})p^{2n}/y }^\infty \frac{dx}{x^{\frac{89}{64}-\varepsilon} }\ll_{N(\mathfrak{N})} \frac{p^{(\frac{23}{32}+\varepsilon)n}}{y^{\frac{39}{64}+\varepsilon}}, \\
II&\ll_{N(\mathfrak{N})}  \sum_{0<m<N(\mathfrak{N})p^{2n}/y} \frac{m^{-\frac{25}{64}+\varepsilon}}{p^{\frac{n}{2}}}\leq \frac{1}{p^{\frac{n}{2}}}\int_0^{ N(\mathfrak{N})p^{2n}/y } \frac{dx}{x^{\frac{25}{64}-\varepsilon}} \ll_{N(\mathfrak{N})} \frac{p^{(\frac{23}{32}+\varepsilon)n}}{y^{\frac{39}{64}+\varepsilon}}.
\end{align*} 
So we are done. 
\end{proof}

We obtain the following estimation on the average:

\begin{proposition}\label{main:thm:1} 
Let $\mathfrak{r}$ be an integral ideal of $F$ coprime to $\mathfrak{p}$.
Then,
\begin{align*}
\lim_{n\rightarrow\infty}L_{\mathrm{av}}(f,\phi,\mathfrak{r})=a_f(\mathfrak{r}).
\end{align*}
\end{proposition}

\begin{proof}
Put $y=p^{\frac{55}{39}n}$ in the estimations of Proposition \ref{1st:gal:av:Lvalue:error:1} and \ref{1st:gal:av:Lvalue:error:2}, then we obtain the desired formula.
\end{proof}

\section{Main results}\label{mainresult:section}

In this section, we obtain the main results.
Let $f\in S_{(\mathbf{k},\mathbf{m}),J}(\mathfrak{N},\chi)$ be a Hecke eigenform. For $\phi\in\Xi_\mathfrak{p}$, let us define an algebraic $L$-value as follows:
$$
L_f(\phi):=\frac{G(\phi)L(\frac{k}{2},f\otimes\phi)}{\Omega_f}\in\mathbb{Q}_f(\phi),
$$
where $\Omega_f$ is the non-zero complex number introduced in Hida \cite{Hida:1994}, which depends on $f$ and the sign of $\phi$, which is totally even.
Let $\mathbb{Q}_f$ be the Hecke field of $f$ over $\mathbb{Q}$, which is a totally real number field or a CM number field due to the rationality of the Hecke eigensystem of $f$ (see Hida \cite[Section 1]{Hida:1994}).
Then, $L_f(\phi)\in\mathbb{Q}_f(\phi)$ and $L_f(\phi)^\sigma=L_{f^\sigma}(\phi^\sigma)$ for any $\sigma\in\mathrm{Gal}(\overline{\mathbb{Q}}/\mathbb{Q})$ by Hida \cite[Theorem 8.1]{Hida:1994}.

Let us recall some important notations and settings:
$\mathfrak{p}$ is a prime ideal of $F$ whose residue degree and ramification degree are $1$, $\mathrm{rank}_{\mathbb{Z}_p}\big(\mathrm{Cl}_F(\mathfrak{p}^\infty)\big)=1$, and $p=N(\mathfrak{p})$ is odd and coprime to $h_F$.
For $\phi\in\Xi_\mathfrak{p}$, let us denote by $\mathfrak{e}(\phi)$ the exponent of the conductor of $\phi$. $n_0>0$ is an integer.
For $\phi\in\Xi_\mathfrak{p}$ with $\mathfrak{e}(\phi)>n_0$, $c_\phi$ is an element of $O_\mathfrak{p}^\times$ such that $\phi(1+\varpi_\mathfrak{p}^{\mathfrak{e}(\phi)-n_0}b)=\mathbf{e}_F(b c_\phi/\varpi_\mathfrak{p}^{n_0})$ for $b\in O_\mathfrak{p}$
and
$$
L_{\mathrm{av}}(f,\phi,\mathfrak{r})=\frac{1}{|G|p^{\mathfrak{e}(\phi)}}\sum_{a\equiv c_\phi(\mathfrak{p}^{n_0})}\mathbf{e}_F\Big(\frac{a}{\varpi_\mathfrak{p}^{\mathfrak{e}(\phi)}}\Big)\sum_{\sigma\in G}\overline{\phi}^\sigma(a\mathfrak{r})\Omega_f L_f(\phi^\sigma),
$$
where $G=\mathrm{Gal}\big(\mathbb{Q}(\phi)/\mathbb{Q}(\mu_{p^{n_0}})\big)$. 
$b_j$ is the element of $1+\mathfrak{p}O_\mathfrak{p}$ such that $b_j^{h_F(p-1)}=\alpha_j$, where $\alpha_j\in F^\times\cap(1+\mathfrak{p}O_\mathfrak{p})$ such that $\mathfrak{a}_j^{h_F(p-1)}=\alpha_j O_F$.

For $S$ an infinite subset of $\Xi_\mathfrak{p}$, there is a sequence $(\phi_m)_m\subset S$ such that $\mathfrak{e}(\phi_m)\rightarrow\infty$ as $m\rightarrow\infty$ due to our assumption.
From now on, for a function $L(\phi)$ of $\phi\in\Xi_\mathfrak{p}$, let us denote by $\lim_{\phi\in S}L(\phi):=\lim_{m\rightarrow\infty}L(\phi_m)$.

\subsection{Determination of newforms by its {\it L}-values}

Firstly, we reprove the non-vanishing result of Kwon-Sun \cite{kwon2020non}:
\begin{theorem}[Non-vanishing of $L$-values]\label{Galois:average:L:value:prop}
$L_f(\phi)\neq 0$ for almost all $\phi\in\Xi_\mathfrak{p}$.
\end{theorem}

\begin{proof}
Clearly, $\mathbb{Q}(\mu_{p^\infty})/\mathbb{Q}_f(\mu_p)\cap\mathbb{Q}(\mu_{p^\infty})$ is infinite Galois and the Galois group is pro-$p$, so that there exists the integer $n_0>0$ such that $\mathbb{Q}(\mu_{p^{n_0}})=\mathbb{Q}_f(\mu_p)\cap\mathbb{Q}(\mu_{p^\infty})$.
Let us assume the contrary, i.e., there exists an infinite subset $S\subset\Xi_\mathfrak{p}$ such that $L_f(\phi)=0$ for each $\phi\in S$. Thus, for $\phi\in S$, $L_f(\phi^\sigma)=0$ for any $\sigma\in G$. 
Therefore, from Proposition \ref{gal:av:nonzero} and \ref{main:thm:1}, we obtain that
$$
0=\lim_{\phi\in S}L_{\mathrm{av}}(f,\phi,O_F)=a_f(O_F),
$$
which is a contradiction as $f$ is an eigenform. So we are done.
\end{proof}

We can determine cusp forms by its $L$-values:
\begin{theorem}[Determination of newforms]\label{determination:Lvalue}
Let $f_i\in S_{(\mathbf{k},\mathbf{m}),J}(\mathfrak{N},\chi_i)$ $(i=1,2)$ be a newform. Let $C$ be a non-zero element of $\mathbb{Q}_{f_1}\mathbb{Q}_{f_2}(\phi)$.  
If $L_{f_1}(\phi)=CL_{f_2}(\phi)$ for infinitely many $\phi\in\Xi_\mathfrak{p}$, then $f_1=f_2$.
\end{theorem}

\begin{proof}
Let us denote $M=\mathbb{Q}_{f_1}\mathbb{Q}_{f_2}(\mu_p,C)$.
Let $n_0>0$ be the integer such that $M\cap\mathbb{Q}(\mu_{p^\infty})=\mathbb{Q}(\mu_{p^{n_0}})$.
By taking the following average
$$
L\mapsto\frac{1}{|G|p^{\mathfrak{e}(\phi)}}\sum_{a\equiv c_\phi(\mathfrak{p}^{n_0})}\mathbf{e}_F\Big(\frac{a}{\varpi_\mathfrak{p}^{\mathfrak{e}(\phi)}}\Big)\mathrm{Tr}_{M(\phi)/M}\big(\overline{\phi}(a\mathfrak{r})L\big)
$$
on the equation $L_{f_1}(\phi)=CL_{f_2}(\phi)$ and taking the limit for $\phi$, we obtain from Proposition \ref{main:thm:1} that
$$
\frac{[\mathbb{Q}(\mu_{p^{n_0}}):\mathbb{Q}_{f_1}(\mu_{p},C)\cap\mathbb{Q}(\mu_{p^\infty})]}{[\mathbb{Q}(\mu_{p^{n_0}}):\mathbb{Q}_{f_2}(\mu_{p},C)\cap\mathbb{Q}(\mu_{p^\infty})]}a_{f_1}(\mathfrak{r})=\frac{C\Omega_{f_1}}{\Omega_{f_2}}a_{f_2}(\mathfrak{r})
$$
for any ideals $\mathfrak{r}$ of $F$ coprime to $\mathfrak{p}$.
Therefore, by the strong multiplicity one, we obtain the desired result.
\end{proof}

\subsection{Generation of cyclotomic Hecke fields by {\it L}-values}
Let us give some useful lemmas to obtain the cyclotomic Hecke field generation:

\begin{lemma}\label{modq:galois:repn} 
Let $f_i\in S_{(\mathbf{k},\mathbf{m}),J}(\mathfrak{N},\chi_i)$ $(i=1,2)$ be a Hecke eigenform.
Let $\varphi_i$ $(i=1,2)$ be a $p$-power order character of modulus $\mathfrak{p}^m$.
Let us assume that $F$ is totally real or CM.
Then, $f_1\otimes\varphi_1\sim f_2\otimes\varphi_2$ implies that $f_1\sim f_2$ and $\varphi_1=\varphi_2$.
\end{lemma}

\begin{proof}
Let us follow the proof of Sun \cite[Proposition 4.2]{sun2018generation}.
Let $\ell\nmid p$ be a rational prime.
Due to our assumption, there exists a Galois representation $\rho_{i}:G_F\rightarrow\mathrm{GL}_2(\overline{\mathbb{Q}}_\ell)$ attached to $f_i$ due to Harris-Lan-Taylor-Thorne \cite{Harris2014on}.
Let $D$ be the compositum of $\overline{F}^{\ker(\rho_1)}$ and $\overline{F}^{\ker(\rho_2)}$, $M$ the ray class field of conductor $\mathfrak{p}^{m}$, and $M_0:=D\cap M$. Then, 
\begin{equation}\label{gal:gp:isom}
\mathrm{Gal}(DM/M_0)\cong\mathrm{Gal}(D/M_0)\times\mathrm{Gal}(M/M_0).
\end{equation}
Let us choose an element $\sigma\in\mathrm{Gal}(DM/M_0)$ such that $\sigma\mapsto(1,\mathfrak{r})$ under the isomorphism (\ref{gal:gp:isom}). 
Since $DM$ is unramified outside $\ell\mathfrak{p}\mathfrak{N}$, there is a prime $v\nmid \ell\mathfrak{p}\mathfrak{N}$ of $F$ such that $\mathrm{Frob}_v|_{DM}$ is an element of the conjugacy class of $\sigma$ in $\mathrm{Gal}(DM/F)$ by Chebotarev's density theorem. 
Then, $\mathrm{Frob}_v$ is trivial in $G_F/\ker(\rho_i)$.
Note that $\mathrm{Gal}(M/M_0)$ is isomorphic to a subgroup of $\mathrm{Cl}_F(\mathfrak{p}^m)$.
Thus, we obtain that
$$
\lambda_1(v)=\lambda_2(v)=2\text{ in }\overline{\mathbb{Q}}_q\text{ and }v=\mathfrak{r}\text{ in } \mathrm{Cl}_F(\mathfrak{p}^{m}),
$$
where $\lambda_i$ is the eigensystem of $f_i$.
Thus, $\varphi_1=\varphi_2$ on $\mathrm{Gal}(M/M_0)$ since $\mathfrak{r}\in\mathrm{Gal}(M/M_0)$ is arbitrary. Note that $M_0$ is an unramified abelian extension of $F$, so that $\mathrm{Gal}(M/M_0)$ is isomorphic to a subgroup of $\mathrm{Cl}_F(\mathfrak{p}^m)$ containing the $p$-part of $\mathrm{Cl}_F(\mathfrak{p}^m)$.
In conclude, $\varphi_1=\varphi_2$, so that $f_1\sim f_2$.
\end{proof}

From now on, we assume that $F$ is totally real or CM.

\begin{lemma}\label{coeff:linear:depend}
Let $m>0$ be an integer and $a$ an element of $O_\mathfrak{p}^\times$.
If
$$
\sum_{\tau\in\mathrm{Hom}_\mathbb{Q}(\mathbb{Q}_f,\mathbb{C})}c_\tau a_{f^\tau}(\mathfrak{r})=0\text{ for any } \mathfrak{r}\in\bigcup_{\kappa\in W_\mathfrak{p} }\bigcup_{i=1}^{h_F}\{r\mathfrak{a}_i^{-1}\mid r\in (a\kappa b_i+\mathfrak{p}^{m}O_\mathfrak{p})\cap\mathfrak{a}_i\},
$$ 
then $c_\tau=0$ for each $\tau\in\mathrm{Hom}_\mathbb{Q}(\mathbb{Q}_f,\mathbb{C})$.
\end{lemma}

\begin{proof}
Let us denote by $\Delta$ the $p$-torsion free subgroup of $\mathrm{Cl}_F(\mathfrak{p}^m)$.
Note that $\mathrm{Cl}_F(\mathfrak{p}^m)\cong\Delta\oplus\mathrm{Cl}_F(\mathfrak{p}^m)_p$, where $\mathrm{Cl}_F(\mathfrak{p}^m)_p$ is the $p$-part of $\mathrm{Cl}_F(\mathfrak{p}^m)$. 
By the orthogonality of characters, we observe from the assumption that
$$
\sum_{\tau\in\mathrm{Hom}_\mathbb{Q}(\mathbb{Q}_f,\mathbb{C})}\sum_{\psi,\varphi}c_\tau\overline{\psi\varphi}(a\kappa b_j\mathfrak{a}_j^{-1})\psi\varphi(\mathfrak{q})a_{f^\tau}(\mathfrak{q})=0
$$
for any prime ideals $\mathfrak{q}\nmid\mathfrak{p}$ of $F$, where $\psi$ and $\varphi$ run through the dual of $\Delta$ and the dual of $\mathrm{Cl}_F(\mathfrak{p}^m)_p$, respectively. This is true for any $\kappa\in W_\mathfrak{p}$ and $j$. 
Let us denote by $a_p$ the image of $a$ in $\mathrm{Cl}_F(\mathfrak{p}^m)_p$.
Then, $aa_p^{-1}\kappa b_j\mathfrak{a}_j^{-1}\in\Delta$ so that 
$\overline{\psi\varphi}(a\kappa b_j\mathfrak{a}_j^{-1})=\overline{\psi}(aa_p^{-1}\kappa b_j\mathfrak{a}_j^{-1})\overline{\varphi}(a_p)$.
Thus, by the orthogonality of characters,
$$
\sum_{\kappa\in W_\mathfrak{p}}\sum_{i=1}^{h_F}\sum_{\varphi}\overline{\psi\varphi}(a\kappa b_i\mathfrak{a}_i^{-1})=
\begin{cases}
h_F(p-1)\sum_\varphi\overline{\varphi}(a_p)&\text{if $\psi=1$ }\\
0&\text{otherwise}
\end{cases},
$$
since $a_p^{-1}a\kappa b_j\mathfrak{a}_j^{-1}$ runs through $\Delta$.
Therefore, we obtain the following from the above equations that 
$$
\sum_{\tau}\sum_{\varphi}c_\tau\overline{\varphi}(a_p)\varphi(\mathfrak{q})a_{f^\tau}(\mathfrak{q})=a_{\sum_{\tau,\varphi}c_\tau\overline{\varphi}(a_p)f^\tau\otimes\varphi}(\mathfrak{q})=0
$$
for any prime ideals $\mathfrak{q}\nmid\mathfrak{p}$ of $F$.
Therefore, we obtain that 
$$
\sum_{\tau,\varphi}c_\tau\overline{\varphi}(a_p)f^\tau\otimes\varphi=0.
$$
Let us assume the contrary, i.e., $c_\tau\neq 0$ for some $\tau$, then $f^{\tau_1}\otimes\varphi_1\sim f^{\tau_2}\otimes\varphi_2$ for some $\tau_1\neq\tau_2$ and $\varphi_1,\varphi_2$ by Lemma \ref{newform:lemma}.
Hence, $f^{\tau_1}\sim f^{\tau_2}$ by Lemma \ref{modq:galois:repn}, which is a contradiction.
\end{proof}

For some $\beta\in\mathbb{Q}_f(\phi)$, the field $\mathbb{Q}(\beta L_f(\phi),\mu_p)$ contains a cyclotomic field:

\begin{proposition}\label{contain:cyclotomic:field}
Let $\beta$ be an element of $\mathbb{Q}_f(\mu_{p^\infty})$. 
Let $n_0>0$ be the integer such that $\mathbb{Q}(\mu_{p^\infty})\cap \mathbb{Q}_f(\beta,\mu_p)=\mathbb{Q}(\mu_{p^{n_0}})$. 
Let us assume that $\mathrm{Tr}_{\mathbb{Q}_f(\mu_{p^{n_0}})/\mathbb{Q}_f}(\beta)\neq 0$. Then,
$\mathbb{Q}(\beta L_f(\phi),\mu_p)\supset \mathbb{Q}(\phi)$ for almost all $\phi\in\Xi_\mathfrak{p}$.
\end{proposition}

\begin{proof}
Let us set $L_\phi:=\mathbb{Q}(\beta L_f(\phi),\mu_p)$.
Let us assume the contrary, i.e., there is an infinite subset $S\subset\Xi_\mathfrak{p}$ such that $L_\phi\not\supset \mathbb{Q}(\phi)$ for each $\phi\in S$.
Note that $S=\bigcup_{c\in (O_\mathfrak{p}/\mathfrak{p}^{n_0}O_\mathfrak{p})^\times}\{\phi\in S\mid c_\phi=c\}$.
Thus, there exists $c\in O_\mathfrak{p}^\times$ such that the set 
$$
S_c:=S\cap\{\phi\in\Xi_\mathfrak{p}:c_\phi=c\}
$$
is infinite.
Let $\mathfrak{r}$ be an integral ideal of $F$ coprime to $\mathfrak{p}$.
Let us investigate the following number: 
$$
A_{a,\mathfrak{r}}(\phi):=\mathrm{Tr}_{\mathbb{Q}_f(\phi)/\mathbb{Q}}\big(\beta\overline{\phi}(a\mathfrak{r})L_f(\phi)\big)
$$
for $\phi\in S_c$.
Note that $\mathbb{Q}_f(\phi)/\mathbb{Q}$ is finite and $\mathbb{Q}_f(\phi)/\mathbb{Q}_f$ is Galois, so that the group $\mathrm{Gal}(\mathbb{Q}_f(\phi)/\mathbb{Q}_f)$ acts on $\mathrm{Hom}_\mathbb{Q}(\mathbb{Q}_f(\phi),\mathbb{C})$ freely. Thus, we have the following decomposition:
$$
\mathrm{Hom}_\mathbb{Q}(\mathbb{Q}_f(\phi),\mathbb{C})=\bigsqcup_{\tau\in\mathrm{Hom}_\mathbb{Q}(\mathbb{Q}_f,\mathbb{C})}\mathrm{Gal}(\mathbb{Q}_f(\phi)/\mathbb{Q}_f)\tau^\prime,
$$
where $\tau^\prime$ is an extension of $\tau$ to $\mathbb{Q}_f(\phi)$, which fixes $\phi$. 
From this, we obtain that
\begin{align*}
A_{a,\mathfrak{r}}(\phi)=\sum_{\tau\in\mathrm{Hom}_\mathbb{Q}(\mathbb{Q}_f,\mathbb{C})}\sum_{\sigma_1\in\mathrm{Gal}(\mathbb{Q}_f(\mu_{p^{n_0}})/\mathbb{Q}_f)}\beta^{\sigma_1\tau^\prime}
\sum_{\sigma\in G}  \overline{\phi}^{\sigma_1\sigma}(a\mathfrak{r}) L_{f^\tau}(\phi^{\sigma_1\sigma}),
\end{align*}
which implies that
\begin{align}\label{average:limit}
\lim_{\phi\in S_c}\frac{1}{|G|p^{\mathfrak{e}(\phi)}}\sum_{a\equiv c(\mathfrak{p}^{n_0})}\mathbf{e}_F\Big(\frac{a}{\varpi_\mathfrak{p}^{\mathfrak{e}(\phi)}}\Big) A_{a,\mathfrak{r}}(\phi)= \sum_{\tau}\frac{\mathrm{Tr}_{\mathbb{Q}_f(\mu_{p^{n_0}})/\mathbb{Q}_f}(\beta)^\tau}{\Omega_{f^\tau}}a_{f^\tau}(\mathfrak{r})
\end{align}
by Proposition \ref{main:thm:1}.
Let $a$ be an element of $O_\mathfrak{p}^\times$ such that $a\equiv c(\mathfrak{p}^{n_0})$.
Then, we obtain that
$$
A_{a,\mathfrak{r}}(\phi)=[\mathbb{Q}_f(\phi):\mathbb{Q}(\phi)L_\phi]\mathrm{Tr}_{L_\phi/\mathbb{Q}}\Big(\beta L_f(\phi)\mathrm{Tr}_{\mathbb{Q}(\phi)/L_\phi\cap \mathbb{Q}(\phi)}(\overline{\phi}(a\mathfrak{r})\big)\Big).
$$
Note that $(\kappa b_j\mathfrak{a}_j^{-1})^{h_F(p-1)}$ is trivial in $\mathrm{Cl}_F(\mathfrak{p}^{\infty})$. Therefore, if $\mathfrak{r}=r\mathfrak{a}_j^{-1}$ for some $r\in c^{-1}(1+p)\kappa b_j+\mathfrak{p}^{n_0}O_\mathfrak{p}$ for $\kappa\in W_\mathfrak{p}$, then $\overline{\phi}(a\mathfrak{r})=\overline{\phi}(1+p)$, which generates $\mathbb{Q}(\phi)$ over $\mathbb{Q}$.
In this case, $A_{a,\mathfrak{r}}=0$ since $L_\phi\cap \mathbb{Q}(\phi)$ is proper in $\mathbb{Q}(\phi)$. Then, we are done by the equation (\ref{average:limit}) and Lemma \ref{coeff:linear:depend}.
\end{proof}

Let us recall that $\mathbb{Q}_f$ is totally real or CM.
We prove the cyclotomic Hecke field generation result for the case of totally real $\mathbb{Q}_f$:

\begin{theorem}[Generation of cyclotomic Hecke fields]\label{themaintheorem}
Let us assume that $f$ is newform, $\mathbb{Q}_f$ is totally real, and $\chi=(\cdot)^{(2-k)\mathbf{t}}$. 
Then, for almost all $\phi\in\Xi_\mathfrak{p}$, there is a character $\phi_0\in\Xi_\mathfrak{p}$ such that
$$
\mathbb{Q}_f(\phi)=\mathbb{Q}\big(\phi_0,A_f(\phi,\phi_0)\big),\text{ where }A_f(\phi,\phi_0):=\frac{L_f(\phi)}{L_f(\phi_0)}\in\mathbb{Q}_f(\phi).
$$
\end{theorem}

\begin{proof}
By Miyake \cite[(1.8)]{miyake1971onautomorphic}, $W_\mathfrak{N}f$ and $f$ have identical eigensystem due to our assumption on $\chi$. Since $W_\mathfrak{N}$ is an involution on the cusp form space 
$S_{(\mathbf{k},\mathbf{m}),J}(\mathfrak{N},\chi)$, 
we can prove that $\pm W_\mathfrak{N}f$ is also a newform. Hence, $W_\mathfrak{N}f=\pm f$ by the strong multiplicity one. Also note that we have the duality (Proposition \ref{Lftn:analytic:conti} and \ref{fricke:twist:commute}) on $L_f(\phi)$, thus we can follow the proof of Sun \cite[Lemma 5.5]{sun2018generation}. Therefore, there exists $\phi_0\in\Xi_\mathfrak{p}$ and an element $\alpha\in\mathbb{Q}(\phi_0)$ such that 
$$
\mathrm{Tr}_{\mathbb{Q}_f(\phi_0)/\mathbb{Q}_f}\Big(\frac{\alpha}{L_f(\phi_0)}\Big)\neq 0.
$$
Note that it is possible to take the reciprocal of the $L$-value by Theorem \ref{Galois:average:L:value:prop}.
Let us denote by $\beta:=\alpha L_f(\phi_0)^{-1}\in\mathbb{Q}_f(\phi_0)$ and $L_\phi:=\mathbb{Q}\big(\phi_0,A_f(\phi,\phi_0)\big)=\mathbb{Q}\big(\phi_0,\beta L_f(\phi)\big)$.
Note that $\mathrm{Tr}_{\mathbb{Q}_f(\phi)/\mathbb{Q}_f}(\beta)\neq 0$ and 
$$
\mathbb{Q}(\phi)\subset L_\phi\subset\mathbb{Q}_f(\phi)
$$ 
for almost all $\phi\in\Xi_\mathfrak{p}$ by Proposition \ref{contain:cyclotomic:field}.
Let $E$ be a Galois number field contains $\mathbb{Q}_f$. For a normal subgroup $H$ of $\mathrm{Gal}(E/\mathbb{Q})$, let us denote by $S_H:=\{\phi\in\Xi_\mathfrak{p}\mid L_\phi\cap E=E^H\}$. Since $E/L_\phi\cap E$ is Galois, the finite collection $\{S_H\}_{H\triangleleft\mathrm{Gal}(E/\mathbb{Q})}$ covers $\Xi_\mathfrak{p}$.
Thus, $|S_H|=\infty$ for some $H\triangleleft\mathrm{Gal}(E/\mathbb{Q})$. 
Therefore, for any $\phi\in S_H$,
$$
L_f(\phi)=\frac{\beta^\sigma}{\beta}L_{f^\sigma}(\phi)\text{ for each } \sigma\in\mathrm{Gal}(E(\phi)/L_\phi).
$$
Note that $f^\sigma$ is also a newform.
From this, we obtain that $f= f^\sigma$ for all $\sigma\in\mathrm{Gal}(E(\phi)/L_\phi)$ by Theorem \ref{determination:Lvalue}, which implies that $\mathbb{Q}_f\subset L_\phi$. So we are done.
\end{proof}

\section*{Acknowledgements}
\thispagestyle{empty}

This work was supported by the Basic Science Research program and the LAMP program of the National Research Foundation of Korea (NRF) grant funded by the Ministry of Education (No. RS-2023-00245291 and RS-2023-00301976). The author would like to thank Dohyeong Kim and Hae-Sang Sun for helpful discussions and comments.

\end{document}